\documentstyle[12pt,leqno]{article}

\newcommand{\finpf}{\ \hfill \rule{2mm}{2mm}\medskip}

\newfam\msbfam
\font\tenmsb=msbm10   \textfont\msbfam=\tenmsb
\font\sevenmsb=msbm7  \scriptfont\msbfam=\sevenmsb
\font\fivemsb=msbm5   \scriptscriptfont\msbfam=\fivemsb
\def\Bbb{\fam\msbfam \tenmsb}

\newcommand{\N}{\mbox{I$\!$N}}

\def\T{{\cal T}}
\def\J{{\cal J}}

\newtheorem{thm}[equation]{Theorem}
\newtheorem{cor}[equation]{Corollary}
\newtheorem{lemma}[equation]{Lemma}
\newtheorem{defn}[equation]{Definition}
\newtheorem{remark}[equation]{Remark}
\newtheorem{pro}[equation]{Proposition}
\newtheorem{Standing Hypotheses}[equation]{Standing Hypotheses}

\date{}
\title{Vector Valued Transference}

\author{E. Berkson, T. A. Gillespie and J. L. Torrea}

\begin{document}
%
%\thanks{The second and third authors were partially supported by  HARP
%network HPRN-CT-2001-00273 of the European Commission.}

\maketitle

{\bf INTRODUCTION}

\

Over the past thirty years, the technique of {\sl transference} has proved extremely effective in obtaining bounds for an operator $T$ acting on a Banach space $X$ (often an $L^p$ space) when $T$ is associated with the representation of a locally compact group $G$ acting on $X$ via a ``convolution-type" formula. The technique has its origins in, {\sl inter alia}, the work of A.P. Calder\'on and A. Zygmund (\cite{CaZy}) on singular integrals and M. Cotlar (\cite{C}) on the ergodic Hilbert transform. An early expository account of the idea of transference was given by Calder\'on (\cite{Ca}) and this was followed by the monograph by R.R. Coifman and G. Weiss (\cite{CW}), who gave a comprehensive survey of the technique and its applications as the theory then stood.

\

The fundamental transference result established by Coifman and Weiss takes the following form. Let $G$ be a locally compact abelian group and $R$ a strongly continuous representation of $G$ in $L^p(\mathcal{M},{\it d} \mu)$ for some measure space $(\mathcal{M}, {\it d} \mu)$ such that there is a uniform bounded $c$ on $\|R_u\| (u\in G)$. In these circumstances, given $k \in L^1(G)$, the formula
$$T_kf = \int_G k(u) R_{-u}fdu $$
defines an operator $T_k$ on $L^p(\mathcal{M}, {\it d} \mu)$, integration being with respect to Haar measure on $G$. The fundamental transference result asserts that 
$\|T_k\| \le c^2 N_p(k)$, where $N_p(k)$ denotes the norm of convolution by $k$ on $L^p(G)$. Roughly speaking, the action of $k$ by convolution on $L^p(G)$ can be transferred to $L^p(\mathcal{M},{\it d} \mu)$ by the representation $R$ of $G$ {\sl with control of norm}. Strictly, Coifman and Weiss consider more generally an amenable group $G$, but most applications involve abelian groups (often $\Bbb{R,T}$ or $\Bbb{Z}$).

\

Since the monograph of Coifman and Weiss appeared, several results involving the transfer of the boundedness of maximal operators and square functions associated to a family of operators have been obtained (see, e.g., \cite{ABG1},\cite{ABG2}) as well as an analogue that applies to representations on an arbitrary Banach space (\cite{BGM}). In this paper we develop a vector valued theory of transference  in which the scalar valued kernel $k$ is replaced by an operator valued kernel. Roughly speaking the principal result is the following (see Theorem  \ref{Principal} below for a precise statement).

\newpage

{\sl Let $X$ and $Y$ be Banach spaces, let $G$ be a locally compact abelian group,
and let $K$ be an operator valued kernel defined on $G$ with values in the space of bounded linear operators from $X$ to $Y$. Suppose that $R$ and $\tilde{R}$ are representations of $G$ on $X$ and $Y$ respectively that intertwine the values of $K$. Then, under suitable boundedness conditions on $R, \tilde{R}$ and $K$, the formula 
$$T_Kx = \int_GK(u)R_{-u}xdu $$
defines a bounded linear operator $T_K$ from $X$ to $Y$ with norm controlled by norm of convolution by $K$ as a mapping from $L^p_X(G)$ into $L^p_Y(G)$, (for all values of $p$ in the range $1\le p < \infty$.)}

\

The layout of the paper is as follows. In section 1, we give a proof of the above vector valued analogue of the Coifman-Weiss transference result; this is in fact a straightforward modification of the proof of the original result. We then show in section 2 how this vector valued result includes the earlier extensions mentioned above. A number of applications to the geometry of Banach spaces are given in section 3. In particular, we give a transference proof of the well known fact that every Banach space in the class {\bf UMD} has non-trivial Rademacher type and cotype (Theorem 3.4) and show that, if a Banach space $X$ has the property that the natural analogue of Rubio de Francia's Littlewood-Paley  inequality for arbitrary intervals (\cite{RdeF2}) is valid for $X$-valued functions, then $X$ has type $p$ for every $p$ in the range $1\le p< 2$. In section 4 several results are proved in the setting of abstract commutative harmonic analysis. In particular, we outline the proof
of the affirmative resolution of a conjecture of Rubio de Francia that was stated in \cite{RdeF1}, the details of which can be found in \cite{BGT}. In the final section of the paper, we indicate how the technique of transference can be used to obtain dimension free estimates for certain operators in an $\Bbb{R}^n$ setting.

\

The notation of the paper is standard and for the most part self-explanatory. We mention only that, given Banach spaces $X$ and $Y$ and a measure space $({\mathcal M},d\mu), {\mathcal L}(X,Y)$ denotes the space of bounded linear operators from $X$ to $Y$ (${\mathcal L}(X)$ when $X=Y$) and, for $1\le p <\infty, 
L^p_X(\mathcal{M})$ denotes the Lebesgue-Bochner space of $p$-integrable $X$-valued functions on $({\mathcal M}, d\mu)$.

\

This survey is the core of a talk given by the third author at the International Conference and 13th Academic Symposium of China on Functional Analysis and Applications. In the treatment below we have, for the first time, presented these results in a complete and unified way. The third author wishes to thank the Department of Mathematics and Statistics of Wuhan University for its warm hospitality, which contributed so much to the success of the Conference

\section {Technical Results}

\begin {thm} \label{Principal} Let $G$ be a locally compact abelian group, let  $X,Y$ be Banach spaces and let $K$ be  a
function in $L^1_{{\cal L}(X,Y)} (G).$ Assume that there exist strongly continuous representations
 $R$ and
$\tilde{R}$ of the group $G$ such that:

(1) for every $u \in G,$ we have $R_u \in {\cal L}(X,X)$ and  $\tilde{R}_u \in {\cal L}(Y,Y);$

(2) there exist constants $c_1$ and $c_2$ such that
$ \|R_u\|_{{\cal L}(X,X)} \le c_1$ and $ \|\tilde{R_u}\|_{{\cal L}(Y,Y)} \le c_2, \,\, u \in G;$

(3) $R$ and $ \tilde{R}$ intertwine  $K,$ in the  sense
$$
K(u)R_v(x) = \tilde{R_v}K(u)(x),  \,\, u,v \in G, \,\, x \in X.
$$

We define the operator $T_K = \int_G K(u)R_{-u}du.$ Then $T_K$ is well defined as an element of
${\cal L}(X,Y)$ and
$$
\|T_K\| \le  \inf _{1 \le p < \infty} (c_1 c_2 N_{p,X,Y}(K)),
$$
where $N_{p,X,Y}(K)$ denotes  the norm of the convolution operator defined by the 
kernel $K$ from $L^p_X$ into $L^p_Y.$
\label{maximal}
\end{thm}

\bigskip

{\bf Proof.} We observe that the operator $T_K$ is well defined as a Bochner integral and
that
\begin{equation}
\label{primera}
\|T_K\| \le  \int \|K(u)\| \|R_{-u}\| du  \le  c_1 \|K\|_{L^1_{{\cal}(X,Y)}}.
\end{equation}

Let $\lambda$  denote Haar measure on $G.$ Assume initially that $K$ has compact support $H \subset G,$
fix $\epsilon > 0$ and let $V$ be an open relatively compact set $V \subset G$ such that
$\frac{ \lambda(V\setminus H)}{ \lambda(V)} < 1+ \epsilon.$  Using (2), (3) and the properties of
the Bochner integral, we have
\begin{eqnarray*}
\|T_Kx\|_Y &\le& c_2 \|\tilde{R}_v \int_G K(u)R_{-u}xdu \|_Y
=  c_2 \| \int_G \tilde{R}_v  K(u)R_{-u}xdu \|_Y  \\
&=& c_2 \| \int_G  K(u)R_{-u}R_vxdu \|_Y
= c_2\| T_K R_vx\|_Y, \,\, x \in X, v \in G.
\end{eqnarray*}
Therefore
\begin{eqnarray*}
\|T_K x\|^p_Y &\le& \frac{c_2^p}{ \lambda(V)} \int_V\| \int_H K(u)R_{v-u}xdu\|^p_Ydv \\
&=& \frac{c_2^p}{ \lambda(V)} \int_V\| \int_G \chi_{V \setminus H}(v-u) K(u)
R_{v-u}xdu\|^p_Ydv  \\
&\le& \frac{c_2^p}{ \lambda(V)} \int_G\| \int_G \chi_{V \setminus H}(v-u)
K(u)R_{v-u}xdu\|^p_Ydv  \\
&\le& \frac{c_2^p}{ \lambda(V)} N_{p,X,Y}(K)^p \int_G\| \chi_{V \setminus H}(w)
R_{w}xdu\|^p_Xdw  \\
&\le& \frac{c_2^p}{ \lambda(V)} N_{p,X,Y}(K)^p \int_{V \setminus H} \|R_wx\|^pdw  \\
&\le& \frac{c_2^p c_1^p N_{p,X,Y}(K)^p}{ \lambda(V)}\lambda(V\setminus H) \|x\|^p_X  \\
&\le& (1+\epsilon) C_1^pc_2^p  N_{p,X,Y}(K)^p \|x\|^p.
\end{eqnarray*}
Letting $\epsilon \rightarrow 0$ gives 
$$\|T_K\|\le c_1c_2N_{p,X,Y}(K)$$
in this case.

Assume now that $ K \in L^1_{{\cal L}(X,Y)}(G).$ Let $K_n$ be  a sequence of compactly
supported functions in $L^1_{{\cal L}(X,Y)}(G)$ such that $\|K_n -K\|_{L^1_{{\cal L}(X,Y)}(G)}
\rightarrow 0.$ For $1\le p <\infty,$
$$N_{p,X,Y}(K_n-K) \le \|K_n-K\|_{L^1_{{\mathcal L}(X,Y)}(G)}$$
and so $\lim_n N_{p,X,Y}(K_n) = N_{p,X,Y}(K)$.  Moreover by (\ref{primera}), $\|T_{K_n}- T_K\| \rightarrow 0$. The desired result now follows by letting $n \rightarrow \infty$ in the inequality 
$$\| T_{K_n} \| \le c_1c_2 N_{p,X,Y}(K_n).$$

\finpf

\bigskip

\section { A unified approach to earlier results}
\setcounter{equation}{0}

In order to see that the vector-valued transference result of the
previous section captures earlier known transference results, it
is necessary to consider various vector-valued extensions of an
operator $S$ defined on a space $L^p(\mu)$ of scalar-valued
functions. To be more precise, let $X$ be a Banach space and
consider when the operator $\tilde{S}=S \otimes Id_X$ defined
initially on the algebraic tensor product $L^p(\mu)\otimes X$ by
\begin{equation} \label{extension} \tilde{S}(\sum \varphi_i
\otimes x_i) = \sum (S\varphi_i) \otimes x_i
\end{equation}
has a bounded extension to the Bochner space $L^p_X(\mu)$. Note that, when $X=\ell^q$, $\tilde{S}$
is given by $\tilde{S}(\{f_j\})=\{Sf_j\}$. Also, when $X=L^p(\nu)$ for some measure $\nu$,
a simple application of Fubini's theorem shows that $\tilde{S}$ does extend to $L^p_X(\mu)$
 with the same norm as $S$. With this observation we obtain from Theorem \ref{maximal} the classical
transference result of Calder\'on and Coifman-Weiss.

\begin{thm} \label{classical} Suppose that $(\Omega,{\cal M},\mu)$ is an arbitrary measure space and let
$u \rightarrow R_u$ be a strongly continuous representation of a locally compact
abelian group $G$ in $L^p(\mu)$, where $1 \le p< \infty$.  Assume that
$c = \sup\{\|R_u\|: u \in G\} < \infty.$
For each $k \in L^1(G)$, use Bochner integration to define the operator $T_k$ of $L^p(\mu)$
into itself by $T_kf = \int_G k(u)R_{-u}fdu.$ Then
$$
\|T_k\| \le c^2 N_{p,L^p(\mu)}(k) =  c^2 N_{p}(k).
$$
\end{thm}
To see this, it is enough to apply Theorem \ref{maximal} with  $X=L^p(\mu)$ and $K=k\otimes Id_X$.
Observing that the convolution operator defined on $L^p(G)$ by
$$
k*g(v) = \int_G k(u)g(v-u)du
$$
has the $L^p(\mu)-$valued extension given by convolution by $K$ and, by the comments above,
this extension has the same norm.

Theorem (\ref{maximal}) also has the following vector valued antecedent (\cite{BG}).

\begin{thm} Suppose that
$u \rightarrow R_u$ is a strongly continuous representation of a locally compact
abelian group $G$ in a Banach space $X$ with
$c = \sup\{\|R_u\|: u \in G\} < \infty.$
For each $k \in L^1(G)$, use Bochner integration to define the operator $T_k$ of $X$
into itself by $T_kf = \int_G k(u)R_{-u}fdu.$ Then, for $1 \le p<\infty$,
$$
\|T_k\| \le c^2 N_{p,X}(k),
$$
where $N_{p,X}(k)$ denotes the norm of convolution by $k $ on $L^p_X(G)$.
\end{thm}

As was observed in \cite{ABG2}, it is not in general possible to transfer strong-type maximal
inequalities under the hypotheses of Theorem \ref{classical}, but such transference is
possible if each operator $R_u$ is assumed to be separation-preserving (that is, $R_uf$ and $R_ug$
have disjoint supports whenever $f$ and $g$ have disjoint supports).

This additional hypothesis can be explained as follows by
considering $\ell^\infty$ extensions.

\begin{lemma} Let $S$ be a bounded linear operator from  $L^p(\mu)$ into itself. Assume
that there exists a bounded positive-preserving (if $f \ge 0 \,\, \mu-$a.e.,
then $ Vf \ge 0 \,\,\mu-$a.e) linear mapping $V$  from  $L^p(\mu)$ into
 itself such that for every $f \in L^p(\mu), |S(f)| \le V(|f|), \mu-$a.e. on $\Omega.$ Then
$S$ has a  $\ell^{\infty}-$valued extension.
\end{lemma}

{\bf Proof.} It is clearly enough to prove that if $V$ is a  bounded positive-preserving
operator then it has a   $\ell^{\infty}-$valued extension. But by using the positive-preserving
property we have
$$
 |V(f_j)(x)| \le V(|f_j|)(x) \le V(\sup_j|f_j|)(x) = V( \|\{f_j\}\|_{l^{\infty}})(x),
$$
then we have
\begin{eqnarray*}
\|\tilde{V}(\{f_j\})\|_{L^p_{l^{\infty}}}
&=& \left( \int (\sup_j|V(f_j)(x)|)^p d\mu(x) \right)^{1/p}
\le \left( \int (V( \|\{f_j\}\|_{l^{\infty}})(x))^p d\mu \right)^{1/p} \\
&\le& \|V\|\left( \int ( \|\{f_j\}\|_{l^{\infty}})(x)^p d\mu \right)^{1/p}
= \|V\| \,\,\|\{f_j\}\|_{L^p_{l^{\infty}}}.
\end{eqnarray*}
\finpf

\begin{cor}\label{exten} Let $S$ be a separation-preserving bounded linear operator from  $L^p(\mu)$ into itself.
Then $S$ has an $\ell^\infty$ extension $\tilde{S}$ and $\|\tilde{S}\|=\|S\|$.
\end{cor}
{\bf Proof.} Since $S$ is separation-preserving, there is a positivity-preserving operator
$|S|$ on  $L^p(\mu)$ such that $|Sf|=|S|(|f|)$ for all $f$ (see \cite{K}). Now apply the above
lemma and note that, by construction, $\|S\|=\|\, | S|\, \| = \|\tilde{S}\| $.
\finpf

\begin{remark} It is worth noting that, conversely, if $S$
is a bounded operator from $L^p(\mu)$ into itself and has a $\ell^\infty-$bounded extension,
then there exists  a bounded positive-preserving linear mapping $V$  from  $L^p(\mu)$ into
 itself such that for every $f \in L^p(\mu), |S(f)| \le V(|f|), \mu-$a.e. on $\Omega.$
See (\cite{V}).
Moreover in this case $S$ has an $X-$bounded extension for any Banach space X. To see this
 we observe that if $f= \sum_{i=1}^n \varphi_i x_i\in X \otimes L^p(\mu)  \otimes X,$ we have
\begin{eqnarray*}
\|\tilde{S}(f) \| &=&\sup_{\xi \in X^*} |<\xi,\tilde{S}(f)>| =
\sup_{\xi \in X^*} |<\xi,\sum_{i=1}^n  S(\varphi_i) x_i>| =
\sup_{\xi \in X^*} |\sum_{i=1}^n  S(\varphi_i)<\xi, x_i> |  \\
&=&
\sup_{\xi \in X^*} |S(<\xi,f>)| \le \sup_{\xi \in X^*} |V(<\xi,f>)|
\le \sup_{\xi \in X^*} V(\|f\|) \le V(\|f\|).
\end{eqnarray*}
\label{extension2}
\end{remark}

\begin{defn} Given $\{k_j\}$ a finite or infinite sequence of functions in $L^1(G).$ For
$1 \le p \le \infty,$ we denote by $N_{p,\infty}(\{k_j\})$ the smallest constant
 $C \in [0, \infty]$ such that
$$
\| \sup_j |k_j*f| \|_{L^p(G)} \le C \|f\|_{L^p(G)},\,\,  \hbox{for all} \,\, f \in L^p(G),
$$
where $*$ denotes convolution with respect to a fixed Haar measure of $G.$
\end{defn}

We are now in a position to establish the following result concerning the transference of
strong-type maximal inequalities (see \cite{ABG2}). (Compare also Theorem (2.11)
of  \cite{BPW} regarding subpositivity and transference of maximal inequalities.)

\begin{thm} Suppose that $R$ is a representation of $G$ in $L^p(\mu)$ such that
$c = \sup\{\|R_u\|: u \in G\} < \infty$ and each $R_u$ is separation-preserving.
Then for  any (finite or infinite) sequence $\{k_j\} \subset
L^1(G)$ such that  $N_{p,\infty}(\{k_j\}) < \infty,$ it follows that the maximal operator $T^*f(x)=
\sup_j|T_{k_j}f(x)|$ is bounded from $L^p(\mu)$ into itself with norm
 \begin{equation}
\| T^* \|_p \le c^2 N_{p,\infty}(\{k_j\}),
\label{A} \end{equation}
where $T_{k_j}$ are the operators defined  by
\begin{equation}
T_{k_j} = \int_G k_j(u)R_{-u}du.
\label{B} \end{equation}
\label{transmax}
\end{thm}

{\bf Proof.} We can assume without loss of generality that the sequence $\{k_j\}$ is
finite. Let $X=L^p(\mu)$ and $Y=L^p_{\ell^\infty}(\mu)$ and define $K(u)  \in {\cal L}(X,Y)$ by
$K(u)f= \{k_j(u)f\}$. We are given the representation $R$ of $G$ in $X$ and, using Corollary
\ref{exten}, we can extend each operator $R_u$ to $Y$ to obtain a representation $\tilde{R}$ of
 $G$ in $Y$. (It is easy to verify that the extended operators do indeed give a representation.)
Notice that, using the notation of Theorem \ref{maximal}, we have $\|T^*\|=\|T_K\|$ and,
 to obtain the conclusion of the present theorem, it suffices to show that
$N_{p,X,Y}(K)=N_{p,\infty}(\{k_j\}).$ This follows from a straightforward application of
 Fubini's theorem.
\finpf

\bigskip

In a similar way, we can consider the transference of square functions from a vector-valued
viewpoint. To do this, recall the Marcinkiewicz-Zygmund theorem which states that, given an
 arbitrary bounded linear operator $S$ from $L^p(\mu)$ into itself (where $1\le p<\infty$),
 the $\ell^2$-valued
extension $\tilde{S}$ of $S$ is also bounded and $\|S\|=\|\tilde{S}\|$.
\begin{defn}Given $\{k_j\}$ a finite or infinite sequence of functions in $L^1(G).$ For
$1 \le p \le \infty,$ we denote by $N_{p,2}(\{k_j\})$ the smallest constant
 $C \in [0, \infty]$ such that
$$
\| (\sum |k_j*f|^2)^{1/2} \|_{L^p(G)} \le C \|f\|_{L^p(G)},\,\,  \hbox{for all} \,\, f \in L^p(G),
$$
where $*$ denotes convolution with respect to a fixed Haar measure of $G.$

\end{defn}

With the Marcinkiewicz-Zygmund theorem in mind, the vector-valued
approach to the transfer of maximal functions can be easily
adapted to give the following theorem \cite{ABG1}.

\begin{thm}Suppose that $R$ is a representation of $G$ in $L^p(\mu)$ such that
$c = \sup\{\|R_u\|: u \in G\} < \infty$ . Then for  any (finite or infinite) sequence
$\{k_j\} \subset L^1(G)$ such that  $N_{p,2}(\{k_j\}) < \infty,$ it follows that the  operator
$Tf(x)=
(\sum_j|T_{k_j}f(x)|^2)^{1/2}$ is bounded from $L^p(\mu)$ into itself with norm
 \begin{equation}
\| T \|_p \le c^2 N_{p,2}(\{k_j\}),
\label{square2} \end{equation}
where $T_{k_j}$ are the operators defined by
\begin{equation}
T_{k_j} = \int_G k_j(u)R_{-u}du.
 \end{equation}
\label{square}\end{thm}

\section {Applications to Geometry of Banach Spaces}
\setcounter{equation}{0}

We recall that a Banach space is said to be of Rademacher type $p$, where $1\le p \le 2,$
(respectively cotype $q$, where $2\le q \le \infty$) if there exists a constant $K$ such that
$$ \int\|\sum_{i=1}^n r_i(t) x_i \|^p dt \le K \sum_{i=1}^n \|x_i\|^p$$
(respectively
$\displaystyle \sum_{i=1}^n \|x_i\|^q \le K \int\|\sum_{i=1}^n r_i(t) x_i \|^q dt $ )
for every finite family of vectors $x_1,\dots,x_n$ in $ X,$  see \cite{LT}. Define the numbers
$$ p(X) = \sup \{p: X \,\hbox{has Rademacher type}\, p \}, \hbox{
and } $$ $$ q(X) = \inf \{q: X \,\hbox{has Rademacher cotype} \, q
\}. $$

\begin{defn} Given Banach spaces $X,Y, $  we shall say that $Y$ is finitely representable
 in $X$ if  for
any finite dimensional subspace $Y_0$ of $Y$ and any $\varepsilon
> 0 $ there exits a finite dimensional subspace $X_0$
of $X$ and a isomorphism ${\cal J}: Y_0 \rightarrow X_0$ with
$\|{\cal J}\|\|{\cal J}^{-1}\| \le 1+ \varepsilon.$ \label{fr}
\end{defn}
The following Theorem, due to B. Maurey and G. Pisier, is well known  and  can be viewed as one
of the  main results in the theory of Banach spaces
( see \cite{MP}).
\begin{thm} The Banach spaces $\ell^{p(X)}$
 and $\ell^{q(X)}$
are finitely representable in $X.$\label{MPP}
\end{thm}

 Let $X$ be a Banach space. We say that $X$ has the {\bf UMD}
 property if the Hilbert transform $H,$ defined initially  in
  $
 L^2({\Bbb R}) \otimes X$   as in (\ref{extension}) has a bounded extension to  $L^2_X({\Bbb R}).$
The {\bf UMD} property was introduced by Burkholder when studying
unconditional convergence of vector-valued  martingale transforms,
it is known that a  Banach space $X$ has the {\bf UMD} property if
and only if the Hilbert transform has a bounded extension to
$L^p_X({\Bbb R})$ for any
 $p$ in the range
$1 < p < \infty.$  See \cite{Bk} and \cite{Bou} and the references there.

\begin{lemma} Assume that $X$ and $Y$ are Banach spaces such that  $X$ has the {\bf UMD}
property and $Y$ is finite representable in
$X$, then $Y$ has the
{\bf UMD} property.
\label{U}
\end{lemma}

{\bf Proof.  } We have to prove that there is a constant  $C$
 such that $$ \|H( \sum_{i=1}^n e_i \varphi_i)
\|_{L^2_{Y}} \le C \| \sum_{i=1}^n e_i \varphi_i\|_{L^2_{Y}}, $$
for all $\{e_i\}_{i=1}^n \subset Y$ and $\{ \varphi_i\}_{i=1}^n
\subset L^2, \, \, n \in {\N}.$  Let $ Y_0 = span\{e_1,\dots,e_n\}
\subset Y, $ and let  $\varepsilon$ be any positive number, then
by using Definition \ref{fr} there exist $X_0$ and an isomorphism
$ {\cal J}:Y_0 \rightarrow X_0 $ such that $\|{\cal J}\|\|{\cal
J}^{-1}\| \le 1+ \varepsilon.$ By using the linearity of the
extension $H$ and the fact that
 $X$ has the {\bf UMD} property,  we have
\begin{eqnarray*}
\left\|H(\sum_{i=1}^ne_i\varphi_i)\right\|_{L^2_{Y}({\Bbb R})} &=& \left\|\sum_{i=1}^ne_i
H (\varphi_i)\right\|_{L^2_{Y}({\Bbb R})} =
\left\|\sum_{i=1}^ne_iH(\varphi_i)\right\|_{L^2_{Y_0}({\Bbb R})}   =\left\|{\cal
J}^{-1}{\cal J}(\sum_{i=1}^ne_iH(\varphi_i))\right\|_{L^2_{Y_0}({\Bbb R})} \\&\le&
\|{\cal J}^{-1}\|\left\|{\cal
J}(\sum_{i=1}^ne_iH(\varphi_i))\right\|_{L^2_{X_0}({\Bbb R})}  = \|{\cal
J}^{-1}\|\left\|\sum_{i=1}^n{\cal J}(e_i)H(\varphi_i)\right\|_{L^2_{X_0}({\Bbb R})}\\&
=& \|{\cal J}^{-1}\|\left\|H(\sum_{i=1}^n{\cal
J}(e_i)\varphi_i)\right\|_{L^2_{X}({\Bbb R})}  \le C_X \|{\cal
J}^{-1}\|\left\|\sum_{i=1}^n{\cal J}(e_i)\varphi_i\right\|_{L^2_{X}({\Bbb R})}\\& =&
C_X \|{\cal J}^{-1}\|\,\left\|{\cal
J}(\sum_{i=1}^ne_i\varphi_i)\right\|_{L^2_{X}({\Bbb R})}  \le C_X \|{\cal
J}^{-1}\|\,\|{\cal J}\|\left\|\sum_{i=1}^ne_i\varphi_i\right\|_{L^2_{Y}({\Bbb R})} \\
&\le& C_X (1+
\varepsilon)\left\|\sum_{i=1}^ne_i\varphi_i\right\|_{L^2_{Y}({\Bbb R})}.
\end{eqnarray*}
\finpf

Now we give a proof of the following well known result  as an
illustration of the use of transference theory in the context of
geometry of Banach spaces.

\begin{thm}
Assume that  a Banach space $X$  satisfies the {\bf UMD}
property. Then $p(X) >1$ and $q(X) < \infty.$
\label{introduction}\end{thm}

{\bf Proof.   } We shall give the proof in the case of $p(X),$ the $q(X)$ case is similar.
 By
Theorem \ref{MPP} and
 Lemma \ref{U}, we have that $\ell^{p(X)}$  satisfies the {\bf UMD}
property. This guarantees that the Hilbert transform maps
$L^2_{\ell^{p(X)}}({{\Bbb R}})$ into itself.  By standard measure theory techniques
this implies that the Hilbert transform $H$ maps
$L^2_{L^{p(X)}({\Bbb R})}({{\Bbb R}})$ into itself.
But this is equivalent to having a uniform bound for the norms on 
$L^2({\Bbb R})_{L^{p(X)}({\Bbb R})}$
  of the operators
 \begin{eqnarray*}
 H_\varepsilon f(x) &=&
  \int_{\{\varepsilon<|x-y|<1/\varepsilon\}}\frac{f(y)}{x-y}dy
  =  \int_{{\Bbb R}}\frac{f(y)}{x-y}\chi_{\{\varepsilon<|x-y|<1/\varepsilon\}}dy\\ &=&\int_{{\Bbb R}}\frac{f(x-y)}{y}\chi_{\{\varepsilon<|y|<1/\varepsilon\}}dy
  = K_\varepsilon*f(x)
  \end{eqnarray*}
  where $$ K_\varepsilon(x) = \frac1{x}\chi_{\{\varepsilon<|x|<1/\varepsilon\}}.$$
Consider the uniformly bounded representation of ${{\Bbb R}}$
 into  $L^{p(X)}({\Bbb R})$ given by $R_u(\varphi)(x) = \varphi(x-u), \quad u\in {\Bbb R}.$
We are in the hypothesis of Theorem \ref{maximal}
with the group $G={{\Bbb R}},$ the Banach spaces  $X=Y=L^{p(X)}({{\Bbb R}})$ and the operators
$K_\varepsilon.$ Therefore the operators $$T_{
K_\varepsilon}\varphi(x)= \int_{{\Bbb R}}
K_\varepsilon(u)R_{-u}\varphi(x)du= K_{\varepsilon}*\varphi(x)$$
are uniformly  bounded in $L^{p(X)}({\Bbb R})$ (with norms controlled
by the norms of $K_\varepsilon$ as operators from
$L^2_{L^{p(X)}({\Bbb R})}({\Bbb R})$ into itself.) But this  implies that
the Hilbert transform is bounded on $L^{p(X)}({{\Bbb R}})$ and hence
that $p(X)>1.$ \finpf

\vspace{0.5 cm}

One of the most remarkable results in recent Harmonic Analysis is the following result due to
Rubio de Francia, see (\cite{RdeF2})

\begin{thm}
 Given an interval $I\subset {\Bbb R},$ we denote by $S_I$ the partial sum operator defined by 
$(S_If)\hat{} = \hat{f} \chi_I,$   where $\hat{f}$ stands for the Fourier transform of the function $f.$
 For every $p$ in the range $ 2 \le p < \infty,$ there exists $C_p >0$
such that, for every sequence $ \{I_k\}$ of disjoint intervals, we have
\begin{equation}
\left\| (\sum_{k} |S_{I_k}f|^2)^{1/2} \right\|_p \le C_p \|f\|_p, \quad  f\in L^p({\Bbb R}).
\label{LP1}\end{equation}
\end{thm}

It is clear that, in order to prove such an inequality, it is enough to prove that there exists
a constant $C_p$ such that, for any finite subfamily $\{I_j\}_{j \in J}$ of ${\cal F},$ we have
\begin{equation}
\left\| (\sum_{j \in J} |S_{I_j}f|^2)^{1/2} \right\|_{L^p({\Bbb R})} \le C_p \|f\|_p, \quad 2 \le p < \infty.
\label{LP2}
\end{equation}

By Kintchine's inequalities, see(\cite{LT}), this last result is equivalent to the existence of
a constant $C_p$ such that for any disjoint family $\{I_j\}_{j \in J}$ and any finite collection
$\{r_j\}_{j \in J}$  of Rademacher functions, we have
\begin{equation}
\| \sum_{j \in J} r_jS_{I_j}f \|_{L^p_{L^p([0,1])}({\Bbb R})} \le C_p \|f\|_p, \quad  2 \le p < \infty.
\label{LP3}\end{equation}
This inequality drives us of the following definition

\begin{defn} Let $X$ be a Banach space and let $2 \le p < \infty$.
We say that $X$ satisfies the
 $LPR_p$
 property if there exists a constant $C_{p,X}$
such that for any finite disjoint family of intervals $\{I_j\}_{j \in J} \subset {\Bbb R}, $  we have:
\begin{equation}
\| \sum_{j \in J} r_jS_{I_j}f \|_{L^p_{L^p_X|0,1]}({\Bbb R})} \le C_{p,X} \|f\|_{L^p_X}({\Bbb R}).
\label{LP4}\end{equation}
\end{defn}

\begin{remark} It is clear that if the Banach space $X$ has the $LPR_p$ property for some $p$ in the range $
2 \le p < \infty,$ then the operators $S_I$ are  bounded from
$L^p_X$ into itself and therefore $X$ must be UMD, see
(\cite{Bou}). An obvious use of Fubini's Theorem says  that the
$L^p$ spaces have the $ LPR_p $ property for $2 \le p < \infty.$
\label{LP5}
\end{remark}

The definition of the $LPR_p$ property has particular significance
for those Banach spaces in which one can defined a notion
 of ``modulus"; that is, for  Banach
lattices. To be precise  we give  the following definition, which   can be found in \cite[p.1 ]{LT}.
\begin{defn} A partially ordered Banach space $E$ over the reals is called a Banach lattice provided:
(i) $x\le y$ implies $x+z\le y+z,$ for every $x,y,z \in E.$ (ii)
$ax\ge 0,$ for every $x\ge 0 $ in $E$ and every nonnegative real
$a.$ (iii) for all $x,y \in E$ there exists a least upper bound
$x\vee y$ and a greatest lower bound $x\wedge y.$ (iv) $\|x\|\le
\|y\|$ whenever $ |x| \le |y| $ where the absolute value $|x|$ of
$x \in E$ is defined by $|x|= x \vee (-x).$
\end{defn}

\begin{defn} A Banach lattice  $E$ is said to be $p$-convex if for every choice of vectors
$\{x_i\}_{i=1}^n$ in $E$ there exists a constant $M$ so that
\begin{itemize}
\item[] $\left\|(\sum_{i=1}^n|x_i|^p)^{1/p}\right\|\le M\left(\sum_{i=1}^n\|x_i\|^p\right)^{1/p},$ \, if
$1 \le p < \infty, \quad$  or
\item[] $\left\|\quad \sup_i|x_i|\quad \right\|\le M \max_i \|x_i\|$ \, if
$ p = \infty. $
\end{itemize}
Analogously a Banach lattice is said to be $q$-concave if there exists a constant $M$ so that
\begin{itemize}
\item[] $\left(\sum_{i=1}^n\|x_i\|^p\right)^{1/p} \le M \left\|(\sum_{i=1}^n|x_i|^p)^{1/p}\right\|\le
$ \, if
$1 \le p < \infty, \quad$  or
\item[]$  \max_i \|x_i\|\le M\left\|\sup_i|x_i|\right\|$ \, if
$ p = \infty. $
\end{itemize}
see \cite[p 46]{LT}
\end{defn}

\begin{remark} The convexity properties of a Banach lattice are closely related to the Rademacher type
and cotype (see \cite[p. 100]{LT}). In particular it is known that if $p(E)>1$ then the lattice is
$p$-convex and $q$-concave for some $1<p<q<\infty$ (see \cite[Corollary 1.f.9]{LT}). Therefore by using
Theorem \ref{introduction} we conclude that if a Banach lattice $E$ is {\bf UMD} then $E$ is
$p$-convex and $q$-concave for some $1<p<q<\infty.$
\label{umdconvexidad}
\end{remark}

Now we present  the following description of the $LPR_p$-property in the case of
Banach lattices.

\begin{lemma}\label{LP6}
 Assume that $E$ is a Banach lattice and $2 \le p < \infty.$ Then the following are
equivalent
\begin{itemize}
\item[(i)] $E$ satisfies the  $LPR_p$-property.
\item[(ii)] There exists a constant $C_p$ such that,  for any disjoint family $\{I_j\}_{j \in J}$ of
intervals in $\Bbb{R}$, we have
\begin{equation}
\label{LP5'}
\left\| (\sum_{j \in J} |S_{I_j}f|^2)^{1/2} \right\|_{L^p_E({\Bbb R})} \le
 C_p \|f\|_{L^p_E({\Bbb R})},
\end{equation}
where $|.|$ is the absolute value in the lattice.
\end{itemize}
\end{lemma}

{\bf Proof.} Observe that both conditions (i) and (ii) imply the  boundedness
of the operators $S_I$ from $L^p_E$ into itself and therefore $E$ must be {\bf UMD}
In order to finish the proof, we shall show that, if a Banach lattice $E$ is {\bf UMD}, then
\begin{equation}
\| \sum_{j \in J} r_jf_j \|_{L^p({\Bbb R})_{L^p_E|0,1]}} \sim
\| (\sum_{j \in J} |f_j|^2)^{1/2} \|_{L^p_E({\Bbb R})}, 
\label{LP7}
\end{equation}
where $1 < p< \infty$ and $\{f_j\}$ is any finite sequence in  $L^p_E.$

To prove (\ref{LP7}), note first that, as  $E$ is {\bf UMD} (see Remark \ref{umdconvexidad}), it must be $q_0$-concave for some $q_0 < \infty.$
 Therefore, if  $q = \max(q_0,p), $ by Jensen's
inequality and $q-$concavity we have
\begin{eqnarray*}
\|\sum_{j \in J}r_jf_j\|^p_{L^p_{L^p_E}({\Bbb R})} &=&
\int_{{\Bbb R}} \int_{|0,1]}  \|\sum_{j \in J}r_j(t)f_j(x)\|_E^pdtdx
\le \int_{{\Bbb R}} (\int_{|0,1]}  \|\sum_{j \in J}r_j(t)f_j(x)\|_E^qdt)^{p/q}dx  \\
&\le& M \int_{{\Bbb R}} \|(\int_{|0,1]}  |\sum_{j \in J}r_j(t)f_j(x)|^qdt)^{1/q}\|_E^pdx
\le M \int_{{\Bbb R}} \|(\sum_{j \in J} |f_j(x)|^2)^{1/2}\|_E^pdx,
\end{eqnarray*}
where, in the last inequality, we have used Kintchine's inequalities.

For the converse, we observe that, as $E$ is {\bf UMD}, then it must be $r_0-$convex for some $r_0 > 1. $
Therefore,  if we put $r= \min(r_0,p), $ by using Kintchine's inequality, the $r-$convexity of $E$ and Jensen's
inequality we have
\begin{eqnarray*}
\int_{{\Bbb R}} \|(\sum_{j \in J} |f_j(x)|^2)^{1/2}\|_E^pdx
&\le& C \int_{{\Bbb R}} |(\int_{|0,1]}  \|\sum_{j \in J}r_j(t)f_j(x)|^rdt)^{1/r}\|_E^pdx \\
 &\le& C \int_{{\Bbb R}} (\int_{|0,1]}  \|\sum_{j \in J}r_j(t)f_j(x)\|_E^rdt)^{p/r}dx  \\
&\le& \int_{{\Bbb R}} \int_{|0,1]}  \|\sum_{j \in J}r_j(t)f_j(x)\|_E^pdtdx.
\end{eqnarray*}
 \finpf

Now we state a Proposition  that is for the $LPR_p$ property the   parallel to Proposition \ref{U} for
the {\bf UMD} property.

\begin{pro} Assume that $X,Y$ are Banach spaces such that  $X$ has the $LPR_p$ property and $Y$ is finite representable in
$X$, then $Y$ has the $LPR_p$ property.
\end{pro}

{\bf Proof. } Let  $R_{X,J}$ be the operator
$$
R_{X,J} = \sum_{j \in J} r_j(t) S_{I_j}f(x).
$$
By definition,   $X$ satisfies the $LPR_p$ property if and only if the operators $R_{X,J}$ are
uniformly bounded from $L^p_{X}({\Bbb R})$ into $L^p_{L^p_{X}([0,1])}({\Bbb R}).$ Therefore we
have  to prove
that there exists a constant $C_p$ (independent of $J$) such that
$$
\left\|R_{Y,J}\left( \sum_{i=1}^n e_i \varphi_i\right) \right\|_{L^p_{L^p_{Y}([0,1])}({\Bbb R})} \le
\left\| \sum_{i=1}^n e_i \varphi_i\right\|_{L^p_{Y}({\Bbb R})},
$$
for all $\{e_i\}_{i=1}^n \subset Y$, $\{ \varphi_i\}_{i=1}^n \subset L^p({\Bbb R})$,  and $n \in {\N}.$ We follow the lines of the proof of Proposition \ref{U}.
Consider $ Y_0 = span \{e_1,\dots,e_n\} \subset Y,$ and $\varepsilon >0$ then by Definition \ref{fr}, there exist a subspace $X_0$ of $X$, and ${\cal J}:Y_0 \rightarrow X_0 $
such that  $\|{\cal J}\|\|{\cal J}^{-1}\| \le 1+ \varepsilon.$ Then, as $X$ has the $ LPR_p$ property,  we
have
\begin{eqnarray*}
\left \|R_{Y,J}(\sum_{i=1}^ne_i\varphi_i)\right\|_{L^p_{L^p_{Y}([0,1])}({\Bbb R})}
&=&\|{\cal J}^{-1}{\cal J}(\sum_{i=1}^ne_iR_J(\varphi_i))\|_{L^p_{L^p_{Y_0}}} \le
\|{\cal J}^{-1}\|\|{\cal J}(\sum_{i=1}^ne_iR_J(\varphi_i))\|_{L^p_{L^p_{X_0}}} \\
&=& \|{\cal J}^{-1}\|\|\sum_{i=1}^n{\cal J}(e_i)R_J(\varphi_i)\|_{L^p_{L^p_{X_0}}} =
\|{\cal J}^{-1}\|\|R_{X,J}(\sum_{i=1}^n{\cal J}(e_i)\varphi_i)\|_{L^p_{L^p_{X}}} \\
&\le& C_p \|{\cal J}^{-1}\|\|\sum_{i=1}^n{\cal J}(e_i)\varphi_i\|_{L^p_{X}}
= C_p \|{\cal J}^{-1}\|\,\|{\cal J}(\sum_{i=1}^ne_i\varphi_i)\|_{L^p_{X}} \\
&\le& C_p \|{\cal J}^{-1}\|\,\|{\cal J}\|\|\sum_{i=1}^ne_i\varphi_i\|_{L^p_{Y}}
\le C_p (1+ \varepsilon)\|\sum_{i=1}^ne_i\varphi_i\|_{L^p_{Y}}.
\end{eqnarray*}
\finpf

Now we give a necessary condition for a Banach space to have the $LPR_p$ property.

\begin{thm}
Assume that $X$ is a Banach space which satisfies the $LPR_p$ property for some
$p,\, 2 \le p < \infty.$ Then $p(X) = 2.$
\label{introduction2}\end{thm}
{\bf Proof.} By Theorem \ref{MPP} and last lemma we have that $\ell^{p(X)}$ satisfies the $LPR_p$
property. Therefore, as we noted in the proof of the last lemma the linear operators
$$
R_{\ell^{p(X)},J}f = \sum_{j \in J} r_jS_{I_j}f
$$
are uniformly bounded from $L^p({\Bbb R})_{\ell^{p(X)}}$ into $L^p({\Bbb R})_{L^p([0,1])_{\ell^{p(X)}}}.
 $ By standard techniques
of extension of operators, see (\cite{RT}), this implies  that the operators
$ R_{L^{p(X)},J}$ are uniformly bounded from  $L^p({\Bbb R})_{L^{p(X)}({\Bbb R})}$ into
 $L^p({\Bbb R})_{L^p([0,1])_{L^{p(X)}({\Bbb R})}}$ and
therefore, as $p(X) \le 2 \le p,$ from $L^p({\Bbb R})_{L^{p(X)}({\Bbb R})}$ into
$L^p({\Bbb R})_{L^{p(X)}([0,1])_{L^{p(X)}({\Bbb R})}}.$ In other words
the operators defined from the space $(L^p \cap L^2)({\Bbb R}) \otimes   L^{p(X)}({\Bbb R})
 $ into
$ (L^p \cap L^2)({\Bbb R})_{L^{p(X)}([0,1])} \otimes L^{p(X)}({\Bbb R})$ by

$$
\sum_ib_i\varphi_i \rightarrow   \sum_ib_i R_J(\varphi_i)  = \sum_ib_i\sum_{j\in J}r_j(.)
(\chi_{I_j} \hat{\varphi_i})\,\check{ },
$$
where $b_i \in L^{p(X)}({{\Bbb R}})$ and $\varphi_i  \in  L^p \cap L^2({{\Bbb R}}),$
have uniform bounded  extensions from $L^p({{\Bbb R}})_{L^{p(X)} ({{\Bbb R}})}$ into
 $L^p({{\Bbb R}})_{L^{p(X)}([0,1])_{L^{p(X)}({{\Bbb R}})}}.
$ Now, by using the ideas in the proof of
 (\cite[Lemma 3.5]{CW}), for each $\chi_{I_j}$ we  find  a sequence of functions $k^n_j$ in
$L^1({{\Bbb R}})$
having compact support such that $\hat{k^n_j}(x) \rightarrow \chi_{I_j}(x),$ a.e., as
 $n \rightarrow \infty. $ We consider  the operators defined from the space
$(L^p \cap L^2)({\Bbb R}) \otimes   L^{p(X)}({\Bbb R})
 $ into
$ (L^p \cap L^2)({\Bbb R})_{L^{p(X)}([0,1])} \otimes L^{p(X)}({\Bbb R})$
by
\begin{eqnarray*}
\sum_ib_i\varphi_i &\rightarrow&
 \sum_ib_i\sum_{j\in J} r_j(.) (\hat{k^n_j}\hat{\varphi_i}) \, \check{ } \,(x) \\
&=&
\sum_ib_i\sum_{j \in J}r_j(.) k^n_j*\varphi_i(x)
= \sum_{jÊ\in J} r_j(.)k^n_j*(\sum_ib_i\varphi_i)(x) \\
&=& (\sum_{jÊ\in J} r_j(.)k^n_j)*(\sum_ib_i\varphi_i)(x)
 = K^n_J * (\sum_ib_i\varphi_i) (x)  \\
&=& \tilde{R}_{L^{p(X)},J}(\sum_ib_i\varphi_i)(x),
\end{eqnarray*}
where $\displaystyle K^n_J \in L^1_{{\cal L}({{\Bbb R}}, L^{p(X)}([0,1]))}$ have compact support and we
are denoting by  $\displaystyle \tilde{R}_{L^{p(X)},J}$ the  convolution operator with $K^n_J.$
By the construction it is easy to see that $\tilde{R}_{L^{p(X)},J}$ have  bounded
extensions from $L^p({{\Bbb R}})_{L^{p(X)}({{\Bbb R}})}$ into
 $L^p({{\Bbb R}})_{L^{p(X)}([0,1])_{L^{p(X)}({ {\Bbb R}})}}.$
We can arrange that the operator norm of $\tilde{R}_{L^{p(X)},J}$ is  bounded by the operator  norm of $ R_{L^{p(X)},J}$  (as in \cite[Lemma 3.5]{CW}). 
Apply  Theorem \ref{maximal} with $R_u\phi(x) = \phi(x-u),
  X =  L^{p(X)}({{\Bbb R}})$
 and $ Y = L^{p(X)}_{ L^{p(X)}({{\Bbb R}})}([0,1]),$ to conclude  that $T_{K^n_j}$ is bounded from $X$ into $Y$
and $\|T_{K^n_j}\| \le \|\tilde{R}_{L^{p(X)},J}\|,$ but
\begin{eqnarray*}
T_{K^n_j}\phi &=& \int {K^n_j}(u)R_{-u}\varphi du = \int \sum_{ j \in J} r_jk^n_j(u)\phi(.-u)du
= \sum_{j \in J}r_j {k^n_j}* \phi.
\end{eqnarray*}
Therefore
$$
 \|\sum_{j \in J}r_j {k^n_j}*\phi\|_{ L^{p(X)}({{\Bbb R}})_{ L^{p(X)}(|0,1])}} \le
\|\tilde{R}_{L^{p(X)},J}\|\|\phi\|_{L^{p(X)}({{\Bbb R}})}
$$
Now, by using the properties of the functions $k^n_j$ and  Fatou's lemma we get
$$
 \|\sum_{j \in J}r_j S_{I_j}\phi\|_{ L^{p(X) ({{\Bbb R}})}_{ L^{p(X)}(|0,1])}} \le
\|\phi\|_{L^{p(X)}({{\Bbb R}})}
$$
 for every finite set $J,$  where $C$ is a bound (independent of $J$) for the norms of    $R_{L^{p(X)},J}$.  This implies that $2 \le p(X)$  and therefore $ p(X) =2.$
\finpf

As we saw in Lemma \ref{LP6},  in certain situations  the existence of a lattice 
structure  in the Banach space can be use to give a
description of a particular geometry property of the space. In
some particular cases the lattice can be described in terms of
classes of functions and then some interesting operators can be
defined. We shall say that  $E$ is a {\bf K\"othe function space}
if $E$ is  a Banach space consisting of
 equivalence classes, modulo equality almost everywhere, of locally integrable real functions on a $\sigma-$finite
 measure space $(\Omega,\Sigma,d\omega),$ such that the following conditions hold:
\begin{itemize}
\item[(i)]
 if $|f(\omega)| \le |g(\omega)|$ a.e. on $\Omega, f$ is measurable and $g \in E,$ then
$f$ belongs to $E$ and $\|f\| \le \|g\|;$
\item[(ii)] for every $A \in \Sigma$ with $\mu(A) < \infty,$ the characteristic function
$\chi_A(\omega)$ of $A$ belongs to $E.$
\end{itemize}  See \cite[p.28]{LT}.

Given a finite subset $J$ of the set ${\Bbb Q}_+$ of the positive rational numbers, we define
$$
{\cal M}_Jf(x) = sup_{r \in J} |B(x,r)|^{-1} \int_{B(x,r)}|f(y)|dy,
$$
where $B(x,r)$ is the ball of radius $r$ centered at $x,$  $ |f(y)| = \sup(f(y),-f(y))$ in the lattice
$E,$ and  $|B(x,r)|$ is the Lebesgue measure of the ball $B(x,r).$
 As an application of our transference Theorem \ref{maximal}
we can prove the following Theorem (see \cite{GMT}).

\begin{thm} The operators ${\cal M}_J$ are not uniformly  bounded from $L^p_{L^1({{\Bbb R}}^n)}({{\Bbb R}}^n)$
into itself for any
$p,  1 < p < \infty.$
\end{thm}
{\bf Proof} Assume that there exists a $p_o$ such that  ${\cal M}_J$ are uniformly bounded from
$L^{p_o}_{L^1({{\Bbb R}}^n)}({{\Bbb R}}^n)$ into itself. This is equivalent to say that the  $\ell^\infty(J)-$
valued operators
\begin{eqnarray*}
N_J f(x) &=& \left\{ |B(x,r)|^{-1} \int_{B(x,r)}|f(y)|dy \right\}_{r \in J}
 = \left\{ |B(0,r)|^{-1} \int_{{{\Bbb R}}^n}| \chi_{B(0,r)}(y)f(x-y)|dy \right\}_{r \in J} \\
&=& \int_{{{\Bbb R}}^n} \left\{ |B(0,r)|^{-1} | \chi_{B(0,r)}(y)\right\}_{r \in J}f(x-y)|dy
= \int_{{{\Bbb R}}^n} \{k_r\}_{r \in J}(y)f(x-y)|dy,
\end{eqnarray*}
are uniformly bounded from $L^{p_o}({{\Bbb R}}^n)_{L^1({{\Bbb R}}^n)}$ into
$L^{p_o}({{\Bbb R}}^n)_{L^1({{\Bbb R}}^n,\ell^\infty(J))}.$
Now we can consider the representation $R$ of ${{\Bbb R}}^n$ into $L^1({{\Bbb R}}^n)$ given by
$R_x(f)(\cdot) =
f(\cdot-x), $ and the representation $\tilde{R}$ of ${{\Bbb R}}^n$ into $L^1_{\ell^\infty(J)}({{\Bbb R}}^n)$
 given by
$\tilde {R}_x(\{f_r\})(\cdot) = \{f_r(\cdot-x)\}_{r \in J}.$ Now we can apply Theorem \ref{maximal} and
 we conclude that the operators
\begin{eqnarray*}
T_{k_J}f(\cdot) = \int_{{{\Bbb R}}^n} \{k_r\}_{r\in J}(u)R_{-u}f(\cdot)du
&=& \left\{\int_{{{\Bbb R}}^n} k_r(u)f(\cdot-u)du \right\}_{r \in J},
\end{eqnarray*}
are uniformly bounded from $L^1({{\Bbb R}}^n)$ into $L^1_{\ell^\infty(J)}({{\Bbb R}}^n). $ This says that the
operators
$$
 M_Jf(x) = sup_{r \in J} |B(x,r)|^{-1} \int_{B(x,r)}|f(y)|dy,
$$ are uniformly bounded from  $L^1({{\Bbb R}}^n)$ into itself,
and therefore we would conclude that the Hardy-Littlewood
 maximal operator is bounded from  $L^1({{\Bbb R}}^n)$ into itself, which is a contradiction.
\finpf

In \cite{GMT} the following definition is given: A {\bf K\"othe
function space} $E$ is said   to satisfy the
 {\bf Hardy-Littllewood property} if there exists
$p_o, 1 < p_o < \infty, $ such that the operators ${\cal M}_J$ are uniformly bounded on
$L^{p_o}_E({{\Bbb R}}^n).$ See \cite{GMT}.

The last  Theorem can be used, together with ideas parallel to the
ideas in Theorems \ref{introduction} and \ref{introduction2},  to prove that if a {\bf K\"othe function space}
 $E$
has the Hardy-Littlewood
property then it must be  $p-$convex , for some $p>1,$ see \cite{GMT}.

\section{Application to functions defined on Groups}
\setcounter{equation}{0}

In this section we give several applications of transference to commutative harmonic analysis, both in a scalar and a vector valued setting.

\begin{thm} Let $G_1$ and $G_2$ be two locally compact abelian groups and let
\newline
 $\pi: G_1 \rightarrow G_2$ be a continuous homomorphism. For $u \in G_1$ and a Banach space $X,$ let $R_u$ denote the isometry from $L^p_X(G_2)$ into $L^p_X(G_2)$ given by
$$
R_uf(t) = f(t + \pi(u)), \,\, ( t \in G_2).
$$
Let $X_1$ and $X_2$ be Banach spaces, let $k \in L^1_{{\cal L}(X_1,X_2)}(G_1),$ and define
$$
T_kf = \int_{G_1} k(u) R_uf du
$$
for  $ f \in L^p_{X_1}(G_2)$. Then $T_k$ is bounded from $L^p_{X_1}(G_2)$ into $L^p_{X_2}(G_2),$ and
$$
\|T_k\| \le N_{p,{X_1},{X_2}}(G_1,k).
$$
where  $N_{p,{X_1},{X_2}}(G_1,k)$ is the operator norm of the  convolution with $k$ from
\newline
 $L^p_{X_1}(G_1)$ into  $L^p_{X_2}(G_1)$
\label{grupos}
\end{thm}

Before proving this theorem, we give the following corollary
\begin{cor} Let $\varphi_1,...,\varphi_N \in L^1({{\Bbb R}}),$ let $p$ and $q$ fixed with 
$1 \le p < \infty ,1 \le q \le \infty ,$  and suppose that
\begin{eqnarray}\label{4.9}
\left\| \, (\sum|\varphi_j*f|^q)^{1/q}\right\|_{L^p({{\Bbb R}})} \le  C \|f\|_{L^p({\Bbb R})} \,\,  \hbox{for all} \,\,
 f \in L^p({{\Bbb R}}).
\label{producto1}
\end{eqnarray}
 Let $S_j$ be the operator on $L^p({\Bbb T})$ corresponding
to the multiplier $\{\hat{\varphi_j}(n)\}_{n \in {\Bbb Z}}.$ Then
\begin{eqnarray}\label{4.8}
\left\| \, (\sum |S_jf|^q)^{1/q}\right\|_{L^p({\Bbb T})} \le  C \|f\|_{L^p({\Bbb T})} \,\,
 \hbox{for all} \,\,
 f \in L^p({\Bbb T}).
\label{producto2}\end{eqnarray}
(When $q = \infty,$ the minorants in the inequalities (\ref{4.9}) and (\ref{4.8}) are interpreted as  $\left\|\sup|\varphi_j*f|\right\|_{L^p({{\Bbb R}})}$  and $\left\|\sup_j|S_jf|\right\|_{L^p({\Bbb T})}$ respectively.)
\end{cor}

{\bf Proof of the Corollary.} Take $G_1 = {{\Bbb R}},\, G_2 = {\Bbb T},$ and
let $ \pi:G_1 \rightarrow G_2$ be given by
$ \pi(u) = e^{iu}.$ Take  $X_1= {\Bbb C}$ and $X_2=\ell^q$ so that ${\cal L}(X_1,X_2) = \ell^q.$ Take
 $k:{{\Bbb R}} \rightarrow \ell^q$ defined by $k(u) = \{ \varphi_j(u)\}_j.$ Then by (\ref{producto1}) we
have $\|k*f\|_{L^p_{\ell^q}({{\Bbb R}})} \le C\|f\|_{L^p_{\Bbb C}}({{\Bbb R}}).$  It follows that the operator
$$
T_kf(e^{it}) = \int_{{\Bbb R}} \{\varphi_j(u)\}_j f(e^{it}te^{iu}) du =
\left\{\int_{{\Bbb R}} \varphi_j(u)\ f(te^{iu}) du \right\}_j
$$
satisfies
$$
\|T_kf\|_{L^p_{\ell^q}({\Bbb T})} \le  C\|f\|_{L^p_{\Bbb C}({\Bbb T})}
$$
In order to prove (\ref{producto2}) we show that $T_kf=\{S_jf\}_j.$ Indeed
\begin{eqnarray*}
\hat{(T_kf)}(n) &=& \left\{ \int_{\Bbb T} \int_{{\Bbb R}} \varphi_j(u)f(e^{it}e^{-iu})e^{-int}dudt \right\}
=\left\{ \int_{\Bbb T} \int_{{\Bbb R}} \varphi_j(u)f(e^{it}e^{-iu})(e^{-iu}e^{it})^{-n}e^{-inu}dudt \right\} \\
&=&\left\{ \int_{\Bbb T} \int_{{\Bbb R}} \varphi_j(u)f(t)e^{-itn}e^{-inu}dudt \right\}
=\left\{ \hat{\varphi_j}(n)\hat{f}(n) \right\} \\
&=& \hat{\{S_jf\}}(n).
\end{eqnarray*}
\finpf

{\bf  Proof of the Theorem (\ref{grupos}).}  We shall apply Theorem (\ref{maximal}) with  $G = G_1,
X =  L^p_{X_1}(G_2), Y = L^p_{X_2}(G_2), R_u = R_u^{X_1}, \tilde{R}_u = R_u^{X_2}$  and 
$K(u) \in L^1_{{\cal L}(L^p_{X_1}(G_2),L^p_{X_2}(G_2))},$ given by
\begin{eqnarray}
\label{Uu}
(K(u)f)(t) = k(u)f(t)  \,\,\, (t \in G_2),
 \end{eqnarray}
in fact
\begin{eqnarray*}
\|K(u)f\|_{L^p_{X_2}(G_2)} &=& \left( \int_{G_2} \|(K(u)f)(t)\|^p_{X_2}dt \right)^{1/p}
= \left( \int_{G_2} \|k(u)f(t)\|^p_{X_2}dt \right)^{1/p}     \\
&\le& \left( \int_{G_2} \|k(u)\|_{{\cal L}(X_1,X_2)} \|f(t)\|^p_{X_1}dt \right)^{1/p} \\
 &=&  \|k(u)\|_{{\cal L}(X_1,X_2)} \|f\|_{L^p_{X_1}(G_2)}.
\end{eqnarray*}
In order to finish the proof we only need to prove that
$$N_{p,L^p_{X_1}(G_2),L^p_{X_2}(G_2)}(G_1,K) =
N_{p,{X_1},{X_2}}(G_1,k),
$$
 but it is clear that any function $F \in L^p(G_1)_{L^p_{X_1}(G_2)}$ can
be realized as a two variables function $F(u,t), u \in G_1,t \in G_2,$ such that for any $u,F(.,t)$
is in $L^p_{X_1}(G_2).$ Then by (\ref{Uu}) we have
\begin{eqnarray*}
\lefteqn{ N^p_{p,L^p_{X_1}(G_2),L^p_{X_2}(G_2)}(G_1,K)
= \sup_{F \in L^p(G_1)_{L^p_{X_1}(G_2)}, \|F\| = 1}
\int_{G_1} \| \int_{G_1} K(u-w)F(w,.)dw\|^p_{L^p_{X_2}(G_2)} du }\\
&=& \sup_{F \in L^p(G_1)_{L^p_{X_1}(G_2)}, \|F\| = 1}
\int_{G_1}  \int_{G_2} \|\int_{G_1} (K(u-w)F(w,.))(t)dw\|^p_{X_2}dt  du \\
&=&\sup_{F \in L^p(G_1)_{L^p_{X_1}(G_2)}, \|F\| = 1}
\int_{G_1}  \int_{G_2}\| \int_{G_1} k(u-w)F(w,t)dw \|^p_{X_2}dt  du \\
 &=&\sup_{F \in L^p(G_1)_{L^p_{X_1}(G_2)}, \|F\| = 1}
\int_{G_2}  \int_{G_1} \|\int_{G_1} \|k(u-w)F(w,t)dw\|^p_{X_2}  du dt  \\
&\le& \sup_{F \in L^p(G_1)_{L^p_{X_1}(G_2)}, \|F\| = 1}
\int_{G_2} \|k*F(.,t)\|^p_{L^p_{X_2}(G_1)}  dt   \\
&\le& \sup_{F \in L^p(G_1)_{L^p_{X_1}(G_2)}, \|F\| = 1}
 N^p_{p,{X_1},{X_2}}(G_1,k) \int_{G_2} \|F(.,t)\|^p_{L^p_{X_1}(G_1)}  dt \\
&=& N^p_{p,{X_1},{X_2}}(G_1,k)
\end{eqnarray*}

For the reverse inequality observe that for any positive $g \in L^p(G_2),$ with
 $\|g\|_{L^p(G_2)}=1,$ we have
\begin{eqnarray*}
N^p_{p,{X_1},{X_2}}(G_1,k)
&=& \sup_{f \in L^p_{X_1}(G_1), \|f\| = 1} \|k*f\|^p_{L^p_{X_2}(G_1)} \\
&=& \sup_{f \in L^p_{X_1}(G_1), \|f\| = 1}
\int_{G_1} \| \int_{G_1} k(u-w)f(w)dw\|^p_{X_2} du \\
&=& \sup_{f \in L^p_{X_1}(G_1), \|f\| = 1}
\int_{G_2} \int_{G_1} \| \int_{G_1} k(u-w)f(w)g(t)dw\|^p_{X_2} dudt \\
&=& \sup_{f \in L^p_{X_1}(G_1), \|f\| = 1}
\int_{G_1} \int_{G_2} \| \int_{G_1} K(u-w)(f(w)g(.))(t)dw\|^p_{X_2} dtdu \\
&\le& \sup_{F \in L^p(G_1)_{L^p_{X_1}(G_2)}, \|F\| = 1}
\int_{G_1} \| \int_{G_1} K(u-w)F(w,.)dw\|^p_{L^p_{X_2}(G_2)} du  \\
&=& N_{p,L^p_{X_1}(G_2),L^p_{X_2}(G_2)}(G_1,K).
\end{eqnarray*}
Where in the penultimate inequality we have used the fact that $F(\omega,t)= f(\omega)g(t)$ belongs to the
space $L^p(G_1)_{L^p_{X_1}(G_2)}$ and $\|F\|_{L^p (G_1)_{L^p_{X_1}(G_2)}}=\|f   \|_{L^p_{X_1}(G_1)} \|g\|_{L^p(G_2)}.$
\finpf

As another application of transference, we discuss briefly the resolution of a conjecture by Rubio de Francia (\cite{RdeF1}). Let $G$ be a compact connected abelian group with dual group $\Gamma$. Then $\Gamma$ can be ordered (in a non-canonical way) so that it becomes an ordered group. Fix any such ordering $\le$ on $\Gamma$ and let $\Gamma^+= \{\gamma\in \Gamma: \gamma \ge 0\}$. A classical result of Bochner \cite{Bo} asserts that, for $1<p<\infty$, the characteristic function of $\Gamma^+$ is a $p$- multiplier. If follows immediately that, for each interval $I$ in $\Gamma$, $\chi_I $ is a $p$-multiplier with a uniform bound on its multiplier norm independent of $I$. (The intervals in $\Gamma$ depend of course, on the particular ordering chosen and may or may not include either of their end-points.)

Given an interval $I$, let $S_I$ denote the corresponding operator on $L^p(G)$. In \cite{RdeF1} J.L. Rubio de Francia observed in the above context that, for $1<p<\infty$ and $\frac1{p}< \frac{2}{q} < \frac{p+1}{p}$, there is a constant $C_{p,q}$ such that
\begin{equation} \label{jl}
\|(\sum_j |S_{I_j} f_j |^q)^{1/q}\|_{L^p(G)} \le C_{p,q} \|(\sum_j | f_j |^q)^{1/q}\|_{L^p(G)}
\end{equation}
for all sequences $\{I_j\}$ of intervals in $\Gamma$ and all sequences $\{f_j\}$ in $L^p(G)$. He noted that, when $G=\Bbb{T}$ or ${\Bbb{T}}^n$, an inequality of the form (\ref{jl}) is in fact valid for all $p,q$ in the range $1<p,q<\infty$ and conjectured 
that this would be the case for an arbitrary compact connected abelian group $G$. this was proved in \cite{BGT} using ideas developed in \cite{BG} and \cite{BGM}, together with Theorem \ref{Principal}. More specifically, a transference argument is used to deduce the result for ${\Bbb{T}}^n$ from that for $\Bbb{T}$. Structural considerations then give the result for a general $G$. This approach has the advantage  of showing that, if $1<p,q<\infty$, then any constant $C_{p,q}$ for which the inequality (\ref{jl}) holds when $G=\Bbb{T}$ will in fact serve fro every $G$ and every ordering in $\Gamma$. In particular, the constant in the inequality for $\Bbb{T}^n$ can be taken to be independent of dimension and of the ordering in ${\Bbb{Z}}^n$. For further details, see \cite{BGT}.

\

We state the result formally as follows.

\begin{thm} Let $1<p,q <\infty$. Then there is a constant $C_{p,q}$ with the following property. For every compact connected abelian group $G$ with ordered dual $(\Gamma, \le )$, 
$$\|(\sum_j |S_{I_j}f_j|^q)^{1/q}\|_{L^p(G)} \le C_{p,q} \|(\sum_j |f_j|^q)^{1/q}\|_{L^p(G)}
$$
for all sequences $\{I_j\}$ of intervals in $\Gamma$ and all sequences $\{f_j\}$ in $L^p(G).$ Furthermore, if $\alpha_{p,q}$ is a constant such that 
$$\|(\sum_j |S_{I_j}f_j|^q)^{1/q}\|_{L^p(\Bbb{T})} \le \alpha_{p,q} \|(\sum_j |f_j|^q)^{1/q}\|_{L^p(\Bbb{T})}
$$
holds for all sequences of  intervals $\{I_j\}$ in $\Bbb{Z}$ and all sequences $\{f_j\}$ in $L^p(\Bbb{T})$, then we can take $C_{p,q}$ to equal $\alpha_{p,q}.$
\end{thm}

The techniques in \cite{BGT}  can also be used to extend Rubio de Francia's Littlewood-Paley inequality discussed in the previous section to the setting of a compact connected abelian group. The precise result is as follows.

\begin{thm} Given $p$in the range $2\le p < \infty,$ there is a constant $C_p$ with the following property. For every compact connected abelian group $G$ with ordered dual $(\Gamma, \le )$,
$$\| (\sum_{I\in {\mathcal F}} |S_If|^2)^{1/2}\|_p \le C_p \|f \|_p, $$
for any family of disjoints intervals  $\mathcal F$ in $\Gamma$. Furthermore the constant $C_p$ can be taken to equal the constant serving for $G=\Bbb{T}.$

\end{thm}

We leave the details to the reader.

\section{Application to dimension free estimates}
\setcounter{equation}{0}

 We consider the operators
$\partial^a(-\Delta)^{-\frac{a+i\gamma}{2}},$ where   $a =0,1 $
 and $\gamma \in {\Bbb R}.$ These operators are  defined for functions whose
Fourier transforms
have compact support by the formula
$$
(\partial^a (-\Delta)^{-\frac{a+i\gamma}{2}}f)\hat{}(\xi) =
 (2\pi i\xi)^a|\xi|^{-(a+i\gamma)}\hat{f}(\xi).
$$
Therefore they have  bounded extensions to $L^2({\Bbb R}^n).$ Since $\Delta$ is
the infinitesimal generator  of the Gauss semigroup, the
operator $(-\Delta)^{-
\frac{a+i\gamma}{2}}$
can also be defined, in terms of the  semigroup, as
$$
(-\Delta)^{-\frac{a+i\gamma}{2}} =
\frac{1}{\Gamma(\frac{a+i\gamma}{2})}\int_0^\infty
e^{t\Delta} \, t^{\frac{a+i\gamma}{2}}\frac{dt}{t},
$$
see \cite{S2}. Therefore, by using the duality in $L^2({\Bbb R}^n),$
the kernels associated, in the sense of Definition \ref{CZ}, with the operators
$\partial^a(-\Delta)^{-\frac{a+i\gamma}{2}}$ as
defined above  can be computed.  In fact, if $f$ is a smooth compactly supported
function, for all $x$ outside the support of $f$ we have

\begin{eqnarray}
\label{gamma}
(-\Delta)^{-\frac{i\gamma}{2}}f(x) &=&
\frac{1}{\Gamma(\frac{i\gamma}{2})}\int_0^\infty \frac{1}{(4\pi t)^{n/2}}
\int_{{\Bbb R}^n}\exp(-\frac{|x-y|^2}{4t})f(y)dy
 \, t^{\frac{i\gamma}{2}}\frac{dt}{t} \\
&=&\frac{\Gamma(\frac{n-i\gamma}{2})}{2^{i\gamma}\pi^{n/2}\Gamma(\frac{i\gamma}{2})}
\int_{{\Bbb R}^n}\frac{1}{|x-y|^{n-i\gamma}}f(y)dy  \nonumber\\
&=&\frac{2^{1-i\gamma}\Gamma(\frac{n-i\gamma}{2})}{\omega_{n-
1}\Gamma(\frac{n}{2})\Gamma(\frac{i\gamma}{2})}
\int_{{\Bbb R}^n}\frac{1}{|x-y|^{n-i\gamma}}f(y)dy \nonumber \\&=&
\Delta_{\gamma}*f(x)
\end{eqnarray}
where $\omega_{n-1} = 2\pi^{n/2}/ \Gamma(\frac{n}{2})$ is the surface area of
the unit sphere
in ${\Bbb R}^n$. Also,

\begin{eqnarray}
\label{nogamma}
\nonumber \partial_{x_j}(-\Delta)^{-\frac{1+i\gamma}{2}}f(x) &=&
\partial_{x_j}\frac{1}{\Gamma(\frac{1+i\gamma}{2})}\int_0^\infty \frac{1}{(4\pi
t)^{n/2}}
\int_{{\Bbb R}^n}\exp(-\frac{|x-y|^2}{4t})f(y)dy
 \, t^{\frac{1+i\gamma}{2}}\frac{dt}{t} \\
&=&-\frac{1}{\Gamma(\frac{1+i\gamma}{2})}\int_0^\infty \frac{1}{(4\pi t)^{n/2}}
\int_{{\Bbb R}^n}\frac{(x_j-y_j)}{2t}\exp(-\frac{|x-y|^2}{4t})f(y)dy
 \, t^{\frac{1+i\gamma}{2}}\frac{dt}{t} \\
&=&-\frac{\Gamma(\frac{n+1-i\gamma}{2})}{2^{i\gamma}\pi^{n/2}\Gamma(\frac{1+i\gamma}{2})}
\int_{{\Bbb R}^n}\frac{2(x_j-y_j)}{|x-y|^{n+1-i\gamma}}f(y)dy \nonumber\\
&=&-\frac{2^{1-i\gamma}\Gamma(\frac{n+1-i\gamma}{2})}{\omega_{n-
1}\Gamma(\frac{n}{2})\Gamma(\frac{1+i\gamma}{2})}
\int_{{\Bbb R}^n}\frac{(x_j-y_j)}{|x-y|^{n+1-i\gamma}}f(y)dy  \nonumber\\&=&
\Delta_{j,\gamma}*f(x). \nonumber
\end{eqnarray}

 Given a kernel $K,$ we denote $K_\varepsilon (x,y) = «
K(x,y)\chi_{\{|x-y|>\varepsilon\}}(x,y).$ The dimension free Theorem is as follows

\begin{thm}
\label{free} Let  $p,\, 1 < p < \infty,$ $\gamma \in {{\Bbb R}},$ and $\alpha$ with $-1 < \alpha < p-1,$
there
exist constant $C_{\alpha,\gamma},$  independent of
$n,$ such that
$$ \left\|\sup_{\varepsilon}|\Delta_{\gamma,\varepsilon} *f| \right\|_{L^p({{\Bbb R}}^n, |x|^\alpha dx)}
 + \left\|\sup_{\varepsilon}\left(\sum_{j}|\Delta_{j,\gamma,\varepsilon} *f|^2\right)^{1/2}\right\|_{L^p({{\Bbb R}}^n,|x|^\alpha dx)}
 \le C_{\alpha,\gamma }
\| f\|_{L^p({{\Bbb R}}^n, |x|^\alpha dx)}. $$
\end{thm}

\vspace{0.5 cm}

In order to prove this theorem, obtained in \cite{GT1}, we shall use some ideas in \cite{DR} \cite{AC}
and \cite{P}. But our intention is to present the proof of the result as an application of a
 weighted transference theory that  can be developed in a  vector valued setting 
(see \cite{GT1} for a  detailed discussion). We  need some
preliminary work to set the stage. As we have indicated earlier, positive
representations constitute an  appropriate vehicle for  the transference of
maximal inequalities. If these representations  satisfy some extra properties,
one can prove some weighted transference Theorems  (see \cite{GT1}
and the references given there). In this section  we  describe briefly such a  weighted transference with vector valued functions and 
present an application. For another application, see \cite{GT2}.

Given a $\sigma-$finite measure space $(X,{\cal F},\mu),$  an
endomorphism of the $\sigma-$algebra ${\cal F}$ modulo null sets
is a set function $\Phi : {\cal F} \rightarrow {\cal F}$ which
satisfies
\begin{enumerate}
\item[(i)] $\Phi ( \cup_n E_n) = \cup_n \Phi( E_n),$ for disjoint $E_n
\in {\cal F}, n = 1,2,\cdots;$
\item[(ii)] $\Phi(X \setminus E) =
\Phi(X) \setminus \Phi(E), $ for all $E \in {\cal F};$
\item[(iii)] given $E \in {\cal F},$ with $\mu(E) = 0,$ then $\mu
(\Phi E) = 0.$ \end{enumerate}
In these circumstances, $\Phi$
induces a unique positive and multiplicative linear operator, also
denoted by $\Phi,$ on the space of (finite-valued or extended)
measurable functions  such that
\begin{equation}\label{measurable}
\Phi(f_n) \rightarrow \Phi(f)\,\, \mu- \hbox{ a.e. whenever  }\, 0
\le f_n \rightarrow f, \, \mu- \hbox{a.e..}
\end{equation}
The action of $\Phi$ on simple functions is given by $$
\Phi(\sum_ic_i\chi_{E_i})(x) = \sum_ic_i\chi_{\Phi(E_i)}(x), \, \,
c_i \in {\Bbb C}. $$
Given a Banach space $B, \Phi$ has an extension, as in (\ref{extension}), to the simple
$B-$valued  functions, also denoted by $\Phi,$ given by $$
\Phi(\sum_i\chi_{E_i}b_i) = \sum_i\chi_{\Phi(E_i)}b_i, \, \, (b_i
\in X, E_i \in {\cal F}). $$ It is clear that, for $f: X
\rightarrow B$ a simple function,
\begin{equation}
\|\Phi( f )(x)\|_B =   \Phi (\|f\|_B)(x), \, ( x \in X). \label{c}
\end{equation}
In other words, if $\Phi$ induces an operator $T$ bounded in
$L^p(\mu),$ then $T$ has a bounded extension, also denoted by $T,$
from $L^p_B(\mu)$ into $L^p_B(\mu)$ for any Banach space $B.$ The
action of $T$ on $L^p(\mu)\otimes B$ is defined as
\begin{equation}T(\sum_i\varphi_i b_i) = \sum_i T(\varphi_i) b_i, \, \, b_i \in B,
\, \, \varphi_i \in L^p(\mu).
\label{vectorextension}
\end{equation}
 The norm of $T$ on $L^p_B(\mu)$
equals the norm of $T$ on $L^p(\mu).$

\begin{Standing Hypotheses}\label{operator}
Throughout, we take $(X,{\cal F},\mu)$ to be a $\sigma-$finite
measure space and $\T = \{T^t: t \in {\Bbb R} \}$ a strongly
continuous one-parameter group  of positive
 invertible linear operators on
$L^p = L^p(X,{\cal F},\mu),$ for some fixed $p$ in the range $1 <
p< \infty,$ such that for each $t \in {\Bbb R},$ there exists a
$\sigma-$endomorphism, $\Phi^t,$ with $T^tf = \Phi^tf.$ In this
case we shall say that $\T$ satisfies ${\bf SH}_p.$
\end{Standing Hypotheses}

 From  the
group structure of $\T,$ it follows that for each $t \in {\Bbb
R},$ there exists  a positive function $J_t$ such that
\begin{equation}
J_{t+s} = J_t \Phi^tJ_s  \hbox{    and    } \int_XJ_t\Phi^tf d\mu
= \int_X f d\mu, \,\, \, \, t,s \in {\Bbb R}. \label{b}
\end{equation}
Using the properties of Bochner integration we have
\begin{equation}
T^t ( \int_{K} T^sfds) =   \int_{K} T^t (T^sf)ds, \, t \in {\Bbb
R} \label{fulanito}
\end{equation}
 for all $f \in L^p(\mu)$ and all compact subsets $K$ of ${\Bbb
 R}.$

\begin{defn} Let $(X,{\cal F},\mu),$  $\T$ and fixed $p$ in the range $1 < p < \infty$ be as
in the ${\bf SH}_p  \, \,5.8$, and let $\omega$ be a
 measurable function on $X$ such that $\omega(x) >0, \,  \, \mu-$almost everywhere.
We shall say that $\omega$ is an {\bf Ergodic $A_p-$weight} with
respect to the group ${\cal T}$ if, for $\mu-$almost all $x \in
X,$ the function $t \rightarrow J_t(x)\Phi^t(\omega)(x)$ is an
$A_p$ weight with an $A_p-$constant independent of $x,$ where $
J_t$ and $\Phi^t$ are as in  (\ref{b}). \label{Ap}
\end{defn}
\noindent We shall denote by  $E_p({\cal T})$ the class of ergodic
$A_p$-weights associated with the group $\T.$
\noindent Given a
weight $\omega$ and a family $\T$ satisfying ${\bf SH}_p\,\, 5.8$, we shall use the following notation
\begin{eqnarray}
\label{notation} \T\omega_x(t) = J_t(x)\Phi^t(\omega)(x).
\end{eqnarray}

 In \cite{GT1} a satisfactory weighted
ergodic theory is developed; one of the outcomes obtained  there is the
following extrapolation result.
\begin{thm}
\label{extrapolation} Let $\T$ be a family of operators satisfying
${\bf SH}_p$ \ref{operator} for every $p$ in the range
$1<p<\infty.$
 Assume that $K$ is a sublinear operator such that
 $\|Kf\|_{L^{p_0}(\omega d\mu)} \le C_\omega\|f\|_{L^{p_0}(\omega d\mu)}$ for every
$\omega \in E_{p_0}(\T),$ where $p_0$ is fixed in the range
$1<p<\infty$ and the constant $C_\omega$ only depends on an
$E_{p_0}(\T)-$constant for $\omega.$ Then $K$ is bounded from
$L^p(\omega d\mu) $ into $L^p(\omega d\mu)$ for every $p, \, 1 <p
< \infty,$ and every $\omega \in E_p(\T).$
\end{thm}

\begin{defn}
Given Banach spaces $B_1, \, B_2,$  and a function
$k \in L^1({{\Bbb R}})_{loc,{\cal L}(B_1,B_2)},$  we shall say that $k$ is a
Calder\'on-Zygmund
kernel if there exists an operator $K$ such that
\begin{itemize}
\item[(i)] for some $p_0, \, 1 < p_0 \le \infty,$  $K$ maps $L^{p_0}_{B_1}({{\Bbb R}})$ into
$L^{p_0}_{B_2}({{\Bbb R}})$;
\item[(ii)] if $\varphi\in L^\infty_{B_1}({{\Bbb R}})$ and has compact support,
then
$$
K\varphi(t) = \int_{{\Bbb R}} k(t-s)\varphi(s)ds, \, \quad t \notin \hbox{support of
} \,  \varphi;
$$.
\item[(iii)] there exists a constant $C$ such that
\begin{itemize}
\item[(iii.1)]$\|k(t)\| \le C|t|^{-1}$ and
\item[(iii.2)] $\|k(t-s) - k(t)\| \le C\frac{|s|}{|t|^2}, \,\, 2|s| < |t|.$
\end{itemize}
\end{itemize}
\label{CZ}
\end{defn}

\begin{remark}
\label{CZ1} Such an operator $K$ is called a Calder\'on-Zygmund operator.
Given $\varepsilon > 0$ we denote by $K_\varepsilon$ the operator obtained by
truncating the kernel
in the standard way, that is
$$
K_\varepsilon \varphi (t) =
\int_{\{s: \varepsilon < |t-s| < 1/\varepsilon\}}k(t-s)\varphi(s)ds.
$$
It is well known that the operator $K^{*}$ defined as $K^{*}\varphi(t) =
\sup_\varepsilon
\|K_\varepsilon \varphi(t)\|_{B_2}$ is bounded
from $L^p_{B_1}({{\Bbb R}}, v) $ into $L^p({{\Bbb R}}, v), \, 1 < p < \infty,$ for
every
$v \in A_p.$ Moreover, the operator norm of $K^{*}$ is majorized by a constant that
 depends only on the operator norm of $K$  on
$L^{p_0}_{B_1}({{\Bbb R}}),$ the constant $C$ in (iii) and  on any  $A_p$ constant
of $v.$
\end{remark}

We now  state the transference  Theorem, whose proof can be found in \cite{GT1}. Recall that $T^t$
has a natural extension to $L^p_{B_1}(X,d\mu),$  also denoted by $T^t$  (see
\ref{extension}).
\begin{thm}
\label{singular integral}
Let $1 < p< \infty$ and let $\T$ be a group of operators satisfying ${\bf
SH}_p \, \,5.8.$ 
Let   $B_1, \, B_2,$ be Banach spaces and $K$ a Calder\'on-Zygmund operator with
 kernel $k$ as in Definition \ref{CZ}.
Given a finite set $\J \subset (0,\infty),$ we define the operator $C_K^\J$ on
$L^p_{B_1}(X,\mu)$ by
$$
C_K^\J f(x) =
\max_{\varepsilon \in \J} \| \int_{\{\varepsilon < |s| < 1/\varepsilon
\}}k(s)T^{-s}f(x)ds\|_{B_2}.
$$
Then
$$
\sup_{\{\J: \J \mbox{ finite } \subset (0, \infty)\}} \|C_K^\J f\|_{L^p(X,\omega)} \le
N_p(K,\T\omega) \|f\|_{L^p_{B_1}(X, \omega)},
$$
for every $\omega \in E_p(\T).$  Here $N_p(K,\T\omega)$  denotes
 an essential bound relative to $x$  of the operator-norm of $K^{*}$ as a
bounded
operator from $L^p_{B_1}({{\Bbb R}},\T \omega_x)$ into
$L^p({{\Bbb R}},\T \omega_x),$ where $\T \omega_x(t)$ is defined in
\ref{notation}.
\end{thm}
We observe that  $\T \omega_x(\cdot) \in A_p$ with an
$A_p$ constant
independent of $x,$ since $\omega \in E_p(\T)$,  and hence such essential bounds exist.

\vspace{0.5 cm}

Let $k$ be  a Calder\'on-Zygmund kernel with the corresponding operator $K$. We consider the
 unit sphere $\Sigma_{n-1}$ of ${{\Bbb R}}^n$ endowed with the rotationally
invariant measure
$d\sigma$ normalized so that $\int_{\Sigma_{n-1}}d\sigma = 1.$
Given a fixed $y' \in  \Sigma_{n-1}$ we consider the one parameter group of
operators
${\T}_{y'} = \{\Phi_{y'}^t\}_t,$ where
$$
\Phi_{y'}^t(f)(x) = f(x+ty'), \, \, x \in {{\Bbb R}}^n, \, t \in {{\Bbb R}}.
$$
Clearly $\|\Phi_{y'}^t(f)\|_{L^p({{\Bbb R}}^n)}= \|f\|_{L^p({{\Bbb R}}^n)}.$
Therefore, if
\begin{equation}
C_{K,\varepsilon,y'} =
\int_{\{\varepsilon < |s| < 1/\varepsilon \}}k(s)\Phi_{y'}^{-s}ds,
\label{preparation}
\end{equation}
then by Theorem \ref{singular integral}
\begin{eqnarray}\label{T}
 \|\{C_{K,\varepsilon,y'}f\}_{\varepsilon \in \J}\|_{L^p_{\ell^\infty(\J)}({{\Bbb R}}^n,\omega)} \le
N_p(K,{\T}_{y'}\omega) \|f\|_{L^p({{\Bbb R}}^n, \omega)}
\end{eqnarray}
for every finite subset $\J$ of $(0,\infty)$ and every  $\omega \in
E_p({\T}_{y'}),$ where
$1 < p< \infty.$

Let $P^0$  and $P^1$ be the projections of the space $L^2(d\sigma)$ into the
subspaces
$H^0$ and $H^1$ of  $L^2(d\sigma)$  generated respectively by the function $1$
and the functions $y'_1,\dots,y'_n.$

\begin{lemma}
\label{hilbert}
 With the notations in \ref{preparation},  we have
$$
 P^0(C_{K,\varepsilon,\cdot}f(x))(y') = C^0_{K,\varepsilon}f(x), \, \, f \in
L^\infty
$$
and

$$
 P^1(C_{K,\varepsilon,\cdot}f(x))(y') = \sum_{j=1}^{n}
C^j_{K,\varepsilon}f(x)Y_j(y'), \,\,  f \in
L^\infty,
$$
where
$$C^0_{K,\varepsilon}f(x) =
\frac{1}{\omega_{n-1}}\int_{\{ z \in {{\Bbb R}}^n : \varepsilon < |z| <
\frac{1}{\varepsilon} \} }
\frac{k(|z|)+k(-|z|)}{|z|^{n-1}}f(x-z)
dz,
$$

$$
C^j_{K,\varepsilon}f(x) = \frac{1}{\omega_{n-1}}\int_{\{ z \in {{\Bbb R}}^n :
\varepsilon < |z| <
\frac{1}{\varepsilon} \} }
\frac{k(|z|)-k(-|z|)}{|z|^{n-1}}f(x-z)
Y_j(\frac{z}{|z|})dz, \, j = 1,\cdots, n
$$
and
$\{Y_j\}_{j=1}^{n}$ are the functions $Y_j(y') = n^{1/2}y_j'$ for $y' \in
\Sigma_{n-1}.$
\end{lemma}

{\bf Proof.  } As $P^1$ is a projection  and $Y_1, \cdots ,Y_n$ are orthonormal
in
$L^2(\Sigma_{n-1},d\sigma),$ we have
$$
 P^1(C_{K,\varepsilon,.}f(x))(y') = \sum_j c_j(x)Y_j(y'),
$$
where
$$
c_j(x) = \int_{\Sigma_{n-1}}C_{K,\varepsilon,y'}f(x)Y_j(y')d\sigma (y').
$$
By using polar coordinates and the fact that the $Y_j's$ are odd functions, the proof
 can be finished.

\finpf

\begin{thm} \label{4.4}Let $K$ be a Calder\'on-Zygmund operator on  ${{\Bbb R}}$
with associated
 kernel $k$ as in \ref{CZ}. Let $ 1 <p < \infty,$ assume that $\omega$ is a
weight
in ${{\Bbb R}}^n$ such that the function $t \rightarrow \Phi^t_{y'}\omega(x)$
is a weight in $A_p({{\Bbb R}})$ with an $A_p$-constant independent of $y'$ and
$x.$
Then there exists a constant $C$ such that
\begin{equation}
\|\, \{\,C^0_{K,\varepsilon}f \, \}_{\varepsilon \in \J}\,
\|_{L^p_{\ell^\infty(\J)}({{\Bbb R}}^n,
\omega)}
\le C\|f\|_{L^p({{\Bbb R}}^n, \omega)}
\end{equation}
and

\begin{equation}
\label{f}
\|\, \{\,(\sum_{j=1}^n|C^j_{K,\varepsilon}f|^2)^{1/2} \, \}_{\varepsilon
\in\J}\,
\|_{L^p_{\ell^\infty(\J)}({{\Bbb R}}^n, \omega)}
\le C\|f\|_{L^p({{\Bbb R}}^n, \omega)}
\end{equation}
for every finite subset $\J$ in $(0,\infty).$   Moreover the constant $C$ can be
taken to be an
upper
bound for the norm of  operators of the form  $K^*: L^p({{\Bbb R}},v) \rightarrow
L^p({{\Bbb R}},v),$
where
$v(t) = \Phi^t_{y'}\omega(x)$ for some $y'$ and $x$ (see Remark \ref{CZ1}).
\end{thm}

{\bf Proof. }
We observe that by using Theorem \ref{extrapolation} it is enough to prove
inequality (\ref{f})
for some $p, \, 1 < p < \infty.$ We shall prove it for $p = 2.$ In fact, using
orthogonality and the representation formula for $P^1$ in  Lemma \ref{hilbert},
we have

\begin{eqnarray*}
\lefteqn{\left\| \, \left\{
\,\left(\sum_{j=1}^n|C^j_{K,\varepsilon}f|^2\right)^{1/2} \,
 \right\}_{\varepsilon \in \J}\,
\right\|_{L^2_{\ell^\infty(\J)}({{\Bbb R}}^n, \omega)} = }\\
&=&  \left\|\, \left\{\,\left(\int_{\Sigma_{n-
1}}|\sum_{j=1}^nC^j_{K,\varepsilon}f
Y_j(y')|^2)d\sigma(y')\right)^{1/2} \, \right\}_{\varepsilon \in \J}\,
\right\|_{L^2_{\ell^\infty(\J)}({{\Bbb R}}^n, \omega)} \\
&=&
\left\|\, \left\{\,\left(\int_{\Sigma_{n-
1}}|P^1(C_{K,\varepsilon,\cdot}f(\cdot))(y') |^2d\sigma
(y')\right)^{1/2}\, \right\}_{\varepsilon \in \J}
\right\|_{L^2_{\ell^\infty(\J)}({{\Bbb R}}^n,
\omega)}\\
&\le&\left\|\, \left\{\,\left(\int_{\Sigma_{n-
1}}|C_{K,\varepsilon,y'}f(\cdot))|^2d\sigma
(y')\right)^{1/2}\, \right\}_{\varepsilon \in \J}
\right\|_{L^2_{\ell^\infty(\J)}({{\Bbb R}}^n,
\omega)}\\
&\le&(\int_{\Sigma_{n-1}} (\|\, \{\,(C_{K,\varepsilon,y'}f(\cdot)\,
\}_{\varepsilon \in \J}
\|_{L^2_{\ell^\infty(\J)}({{\Bbb R}}^n, \omega)} )^2d\sigma (y'))^{1/2}\\
&\le&
(\int_{\Sigma_{n-1}} N_2(K,\T_{y'}\omega)^2 \|\, f \|^2_{L^2({{\Bbb R}}^n,\omega))}d\sigma(y'))^{1/2}
\\
&\le& C \|\, f \|_{L^2({{\Bbb R}}^n,\omega))},
\end{eqnarray*}
where in the penultimate  inequality we have used  \ref{T}. The case $C^0_{K,\varepsilon}$ is
simpler and we live the details to the reader.
\finpf

\begin{cor}
\label{pesos}
Let $  1 <p < \infty$ and let $
-1 < \alpha < p-1.$
Let $K$ be a Calder\'on-Zygmund operator on  ${{\Bbb R}}$ with associated
 kernel $k$ as in \ref{CZ} and consider the operators $C^0_{K,\varepsilon}$ and
$C^j_{K,\varepsilon}$ defined in Lemma \ref{hilbert}. Then there exists a
constant $C_{\alpha,p}$
such that
\begin{equation}
\int_{{{\Bbb R}}^n}
\sup_{\varepsilon>0} |C^0_{K,\varepsilon}f(x)|^p |x|^\alpha dx
\le C_{\alpha,p} \int_{{{\Bbb R}}^n} |f(x)|^p|x|^{\alpha}dx,
\end{equation}
and
\begin{equation}
\label{ff}
\int_{{{\Bbb R}}^n}
\sup_{\varepsilon>0}(\sum_{j=1}^n|C^j_{K,\varepsilon}f(x)|^2)^{p/2} |x|^\alpha
dx
\le C_{\alpha,p} \int_{{{\Bbb R}}^n} |f(x)|^p|x|^{\alpha}dx.
\end{equation}
\end{cor}

{\bf Proof. } In order to use  Theorem \ref{4.4}, it will be enough to show that,
given
$x \in {{\Bbb R}}^n$ and $y' \in \Sigma_{n-1},$ the function $t \rightarrow
|x+ty'|^{\alpha}$
is an $A_p$-weight on ${{\Bbb R}},$ with an $A_p$-constant independent of $x$
and $y'.$

Fix $x \in {{\Bbb R}}^n, y' \in \Sigma_{n-1}$ and decompose $x$ as $x = x_1 +
t_0y',$
with $x_1 \bot y'.$ Then, as $|y'| = 1,$ we have
$|x+ty'| = (|x_1|^2 + |t_0+t|^2)^{1/2} \sim |x_1| + |t_0+t|.$ Therefore
$|x+ty'|^{\alpha} \sim |x_1|^\alpha + |t_0 + t|^\alpha.$ Hence if $M$ is
 the Hardy-Littlewood maximal
operator and we denote by $\varphi_s$ the translate function
$\varphi_s(t) = \varphi(t-s),$ by using the translation properties of Lebesgue
measure
and the fact that $|t|^\alpha$ is a $A_p$-         weight,
we have
\begin{eqnarray*}
\int_{{\Bbb R}}|M\varphi(t)|^p(|x_1|^\alpha + |t_0 + t|^\alpha)dt &=&
\int_{{\Bbb R}}|M\varphi(t)|^p|x_1|^\alpha dt +\int_{{\Bbb R}}|M\varphi(t)|^p|t_0 +
t|^\alpha dt \\
&=& |x_1|^\alpha \int_{{\Bbb R}}|M\varphi(t)|^p dt +\int_{{\Bbb R}}|M\varphi_{t_0}(t)|^p|t|^\alpha dt
\\
&\le& |x_1|^\alpha C_p \int_{{\Bbb R}}|\varphi(t)|^p dt
+A_p(|t|^{\alpha})\int_{{\Bbb R}}|\varphi_{t_0}(t)|^p|t|^\alpha dt \\
&\le& (C_p + A_p(|t|^\alpha)) \int_{{\Bbb R}}|\varphi(t)|^p(|x_1|^\alpha+|t_0+t|^\alpha) dt.
\end{eqnarray*}
It follows that $|x_1|^\alpha +|t_0+t|^\alpha,$ and hence $|x+ty'|^\alpha,$ is
an $A_p$-weight
with an $A_p$-constant on ${{\Bbb R}}$ independent of $x$ and $y'.$
\finpf

{\bf Proof of  Theorem \ref{free}. }
 We consider the Calder\'on-Zygmund operators
on
${{\Bbb R}}$ given by the Calder\'on-Zygmund kernels
$k_0 (t) =
|t|^{-1+i\gamma}$ with $\gamma \ne 0$ and
$k_1(t) =  t^{-1}$   (see \cite[Ch II]{S1}).
Therefore with this notation we have in Lemma \ref{hilbert}

\begin{eqnarray*}
C^0_{K_0,\varepsilon}f(x) &=& \frac{1}{\omega_{n-1}}
\int_{\{z \in {{\Bbb R}}^n: \varepsilon < |z| < \frac{1}{\varepsilon} \}}
\frac{2k_0(|z|)}{|z|^{n-1}}f(x-z)dz \\
&=& \frac{2}{\omega_{n-1}} \int_{\{z \in {{\Bbb R}}^n: \varepsilon < |z| <
\frac{1}{\varepsilon} \}}
\frac{1}{|z|^{n-i\gamma}}f(x-z)dz.
\end{eqnarray*}
In other words, by \ref{gamma}, $C^0_{K_0,\varepsilon}f(x) =\kappa_n \Delta_{\gamma,
\varepsilon}*f(x)$
with $\kappa_n=\frac{2^{i\gamma}\Gamma(\frac{n}{2})\Gamma(\frac{i\gamma}{2})}
{\Gamma(\frac{n-i\gamma}{2})}.$ Using
Stirling's formula it is easy to see that $|\kappa_n| \sim C.$ Therefore
Corollary
\ref{pesos} applies and we obtain
\begin{eqnarray*}
\int_{{{\Bbb R}}^n}
\sup_{\varepsilon> 0}|\Delta_{\gamma, \varepsilon}f(x)|^p |x|^\alpha dx
\le C_{\alpha,p} \int_{{{\Bbb R}}^n} |f(x)|^p|x|^{\alpha}dx.
\end{eqnarray*}

On the other hand, by Lemma \ref{hilbert} we also have
\begin{eqnarray*}
C^j_{K_1,\varepsilon}f(x) &=& \frac{2n^{1/2}}{\omega_{n-1}}
\int_{\{z \in {{\Bbb R}}^n: \varepsilon < |z| < \frac{1}{\varepsilon} \}}
 \frac{1}{|z|}
f(x-z)\frac{z_j}{|z|^n}dz,
\end{eqnarray*}
and so, by (\ref{nogamma}) with $\gamma = 0$,  $C^j_{K_1,\varepsilon}f(x) = -\kappa_n
\Delta_{j,\varepsilon}*f(x),$ where $\kappa_n =
\frac{n^{1/2}\Gamma(\frac{n}{2})\Gamma(\frac{1}{2})} {\Gamma(\frac{n+1}{2})}
\Delta_{j,\varepsilon}.$ As before, Stirling's formula gives $|\kappa_n| \sim C$
and therefore
the case $m=1$ and $\gamma = 0$  in the theorem follows from Corollary
\ref{pesos}.

\

{\bf Acknowledgment.}  The second and third authors were partially supported by  HARP
network HPRN-CT-2001-00273 of the European Commission.

\
EB:Department of Mathematics, University of Illinois, 1409 West Green Street, Urbana, Illinois 61801, USA.\,\,   berkson@illinois.edu

TAG:School of Mathematics, University of Edinburgh, James Clerk Maxwell Building, Edinburgh EH9 3JZ, Scotland. \,\, T.A.Gillespie@ed.ac.uk

JLT: Departamento de Matem\'aticas, Universidad Aut\'onoma de Madrid, 28049 Madrid, Spain. \,\,  joseluis.torrea@uam.es


\begin{thebibliography}{ABC}





\bibitem [A,B,G1]{ABG1} N.Asmar, E. Berkson and T.A. Gillespie, {\sl Transferred bounds
for square functions},
Houston J. Math. (1991), 525-550.





\bibitem [A,B,G2]{ABG2} N.Asmar, E. Berkson and T.A. Gillespie, {\sl Transference of
strong type maximal inequalities by separation-preserving representations},
American J. Math. 113 (1992), 47-74.


\bibitem [A,C]{AC} P. Auscher and M.J. Carro, {\sl Transference for radial
multipliers
and dimension free estimates},
Trans. Amer. Math. Soc. 342 (1994), 575-593.

\bibitem [B,G]{BG}  E. Berkson and  T.A. Gillespie, {\sl The generalized M. Riesz Theorem and transference}, Pacific J. Math. 120 (1985), 279-288.

\bibitem [B,G,M]{BGM}  E. Berkson, T.A. Gillespie and P.S. Muhly {\sl Generalized analyticity in UMD
spaces}, Arkiv f\"or Matematik, 27 (1989),1-14.



\bibitem [B,G,T]{BGT} E. Berkson, A. T. Gillespie, J.L. Torrea, {Proof of a conjecture of
 Jos\'e L. Rubio de Francia ,} Math. Z. 251 (2005), 285-292. 



\bibitem [B,P,W]{BPW} E, Berkson, M. Paluszy\u{n}ski and G. Weiss, {\sl Transference couples and their applications to convolution operators and maximal operators,} Proc. of the Conference on the Interaction between Functional Analysis, Harmonic Analysis and Probability (Univ. of Missouri at Columbia, 1994), Lecture Notes in Pure and Appl. Math. 175, Marcel Dekker, New York (1996), 69-84.

\bibitem [Bo]{Bo} S. Bochner {\sl Additive set functions on groups} Annals of Math. 40
(1939), 769-799


\bibitem [Bou]{Bou} J. Bourgain, {Some remarks on Banach spaces in which martingale sequences are unconditional}, Arkiv f\"or Matematik 21 (1983),163-168.

\bibitem [B]{Bk}  D. Burkholder {\sl A geometric condition that implies the
existence of certain singular integrals of Banach-spaces valued functions},
Conference on Harmonic Analysis in Honor of A. Zygmund 1981, Wodsworth
International Group, Belmont, California 1 (1983), 270-286.


\bibitem [C]{Ca}  A. P. Calder\'on. {\sl Ergodic Theory and translation invarian operators,} Proc. Nat. Acad. Sci.
 (U.S.A. ) 59 (1968), 349-353.
 
 \bibitem [C,Z]{CaZy}  A. P. Calder\'on  and A. Zygmund. {\sl On singular integrals,} Amer. J. Math.
78 (1968), 289-309.





\bibitem [C,W]{CW} R.R. Coifman and G. Weiss. {\sl Transference methods in Analysis} C.B.M.S. Regional Conf. Series
in Math. No.31, Amer. Math. Soc., Providence, R.I., 1977, reprinted 1986

\bibitem [Co] {C} M. Cotlar. {\sl A unified theory of Hilbert transforms and ergodic theorems}
Rev. Mat. Cuyana1 (1955), 105-167.

\bibitem[D,R]{DR} J. Duoandikoetxea and J.L. Rubio de Francia, {\sl Estimations
ind\'ependantes de la
dimension
pour les transform\'ees de Riesz,} C.R. Acad. Sci. Paris S\'er. I 300 (1985),
193-
196.


\bibitem [GC,M,T]{GMT} , J. Garc\'{\i}a-Cuerva, R. Mac\'{\i}as and J.L. Torrea {\sl The
Hardy-Littlewood property of Banach lattices.} Israel J. of Math. 83 (1993),177-201.

\bibitem [G,T1]{GT1} A. T. Gillespie, J.L. Torrea, {\sl Weighted
ergodic theory and dimension free estimates,} 
Q.J. Math.  54 (2003), 257-280.

\bibitem [G,T2]{GT2}  A. T. Gillespie, J.L. Torrea, {\sl Dimension free estimates
for the oscillation of Riesz transforms,} 
Israel.J. Math. 141 (2004),  125–144.


\bibitem [K]{K} C.-H. Kan, {\sl Ergodic properties for Lamperti operators},
Canadian J. Math. 30 (1978), 1206-1214.


\bibitem [L,T]{LT} J. Lindenstrauss and L. Tzafriri {\sl Classical Banach Spaces II.
Function spaces.} Springer-Verlag, Berlin 1979.

\bibitem [M,P]{MP} B. Maurey and G. Pisier ,{\sl Series de variables aleatoires vectorielles independantes et
 proprietes geometriques des spaces de Banach.} Studia Math. 58 (1976), 45-90.








\bibitem [M]{M} P.A. Meyer, {\sl Note sur le processus d'Ornstein-Uhlenbeck},
Springer Lecture
Notes
in Mathematics 920 (1982), 95-132.








\bibitem [P]{P} G. Pisier, {\sl Riesz transforms: a simpler analytic proof of
P.A.
 Meyer's inequality}, Lecture Notes in Math., vol. 1321, Springer-Verlag, 1988,
485-501.




\bibitem [RdeF1]{RdeF1} J.L.Rubio de Francia {\sl Continuity and Pointwise Convergence of Operators
in Vector Valued $L^p$ Spaces}
Bulletino Unione Mat. Italiana 17 (1980), 650-660.

\bibitem [RdeF2]{RdeF2} J.L.Rubio de Francia {\sl A Littlewood-Paley inequality for
 arbitrary intervals,} Rev. Mat.
Iberoamericana 2 (1985), 1-4


\bibitem [RdeF,T]{RT} J.L.Rubio de Francia and J.L Torrea  {\sl Some Banach techniques
in vector-valued Fourier Analysis} Colloquium  Math. 54 (1987), 273-284.

\bibitem[S1]{S1} E.M.Stein , {\sl Singular Integrals and Differentiability
Properties of Functions}, Princeton Univ. Press, Princeton, N.J., 1970.


\bibitem [S2]{S2} E.M. Stein, {\sl Topics in Harmonic Analysis related to the
Littlewood-Paley
Theory},
Annals of Mathematical Studies, Princeton Univ. Press, Princeton, N.J., 1970.


\bibitem[V]{V} B. Virot {\sl Extensions vectorielles d'operateurs lineairs bornes
sur $L^p$.} Seminaire Maurey-Pisier 1980-81. Expos\`e no.7

\end{thebibliography}
\end{document}